\RequirePackage{tikz}

\documentclass[12pt]{amsart}

\usepackage{geometry}
\usepackage{hyperref}

\usepackage[normalem]{ulem}

\usepackage[style=numeric, sorting=nty]{biblatex}

\bibliography{library}

\usepackage{amsfonts}
\usepackage{amsmath}
\usepackage{amssymb}
\usepackage{mathtools}
\usepackage{amsthm}

\newcommand{\sbm}[1]{{\let\amp=&\left[\begin{smallmatrix}#1\end{smallmatrix}\right]}}

\usepackage[all]{xy}
\usetikzlibrary{cd}
\usetikzlibrary {arrows.meta}
\usetikzlibrary{arrows}
\usetikzlibrary{calc}
\tikzset{circle through 3 points/.style n args={3}{%
insert path={let \p1=($(#1)-(#2)$),\p2=($(#1)!0.5!(#2)$),
    \p3=($(#1)-(#3)$),\p4=($(#1)!0.5!(#3)$),\p5=(#1),\n1={(-(\x2*\x3) + \x3*\x4 + \y3*(-\y2 +
    \y4))/(\x3*\y1 - \x1*\y3)},\n2={veclen(\x5-\x2-\n1*\y1,\y5-\y2+\n1*\x1)} in
    ({\x2+\n1*\y1},{\y2-\n1*\x1}) circle (\n2)}
}}
\tikzset{half circle through 3 points/.style n args={3}{%
insert path={let \p1=($(#1)-(#2)$),\p2=($(#1)!0.5!(#2)$),
    \p3=($(#1)-(#3)$),\p4=($(#1)!0.5!(#3)$),\p5=(#1),\n1={(-(\x2*\x3) + \x3*\x4 + \y3*(-\y2 +
    \y4))/(\x3*\y1 - \x1*\y3)},\n2={veclen(\x5-\x2-\n1*\y1,\y5-\y2+\n1*\x1)} in
    ({\x2+\n1*\y1},{\y2-\n1*\x1}) ++(0:\n2) arc(0:180:\n2)}
}}
\tikzset{ray through 2 points/.style n args={2}{%
insert path={let \p1=(#1),\p2=(#2) in
    (\p1) -- ({\x1 + 2 * (\x2 - \x1)},{\y1 + 2 * (\y2 - \y1)})}
}}
\tikzset{long ray through 2 points/.style n args={2}{%
insert path={let \p1=(#1),\p2=(#2) in
    (\p1) -- ({\x1 + 4 * (\x2 - \x1)},{\y1 + 4 * (\y2 - \y1)})}
}}
\usepackage{graphicx}
\usepackage{float}
\usepackage{caption}
\usepackage{subcaption}
\usepackage{xcolor}
\usepackage{colortbl}

\usepackage{enumitem}

\theoremstyle{thmstyleone}%
\newtheorem{theorem}{Theorem}[section]
\newtheorem{lemma}[theorem]{Lemma}%
\newtheorem{claim}[theorem]{Claim}%

\theoremstyle{thmstyletwo}%
\newtheorem{example}[theorem]{Example}%

\theoremstyle{thmstylethree}%
\newtheorem{definition}[theorem]{Definition}%

\newcommand{\fr}[1]{{\color{purple} #1}}

\newcommand{\rA}{\mathrm{A}}

\newcommand{\bZ}{\mathbb{Z}}

\newcommand{\bR}{\mathbb{R}}
\newcommand{\bV}{\mathbb{V}}
\newcommand{\bW}{\mathbb{W}}
\newcommand{\bX}{\mathbb{X}}

\newcommand{\cA}{\mathcal{A}}
\newcommand{\cB}{\mathcal{B}}

\newcommand{\cI}{\mathcal{I}}

\newcommand{\cN}{\mathcal{N}}
\newcommand{\cE}{\mathcal{E}}

\newcommand{\cS}{\mathcal{S}}
\newcommand{\cU}{\mathcal{U}}

\newcommand{\EXT}{\textrm{EXT}}
\newcommand{\rB}{\textrm{B}}
\newcommand{\rD}{\textrm{D}}

\newcommand{\lbirths}{\texttt{lbirths}}
\newcommand{\ldeaths}{\texttt{ldeaths}}

\newcommand{\Id}{\textrm{Id}}

\newcommand{\redCell}{\cellcolor{red!20}}
\newcommand{\blCell}{\cellcolor{blue!20}}

\newcommand{\TX}{\mathrm{T}(\bX)}
\newcommand{\DX}{\mathrm{D}(\bX)}
\newcommand{\DXi}[1]{\mathrm{D}(\bX_{#1})}
\newcommand{\VorX}{\mathcal{V}(\bX)}
\newcommand{\AX}[1]{\mathrm{A}_{#1}(\bX)}
\newcommand{\AXn}{\mathrm{A}(\bX)}

\newcommand{\Ho}{\mathrm{H}}

\newcommand{\PH}{\mathrm{PH}}

\newcommand{\Tot}{\textrm{Tot}}

\newcommand{\bOne}{\textbf{1}}

\DeclareMathOperator{\Ima}{Im}
\DeclareMathOperator{\Ker}{Ker}

\DeclareMathOperator{\PVect}{PVect}
\DeclareMathOperator{\Cycle}{Cycle}

\newcommand{\bigBoxPlus}{
  \mathop{
    \vphantom{\bigoplus}
    \mathchoice
      {\vcenter{\hbox{\resizebox{\widthof{$\displaystyle\bigoplus$}}{!}{$\boxplus$}}}}
      {\vcenter{\hbox{\resizebox{\widthof{$\bigoplus$}}{!}{$\boxplus$}}}}
      {\vcenter{\hbox{\resizebox{\widthof{$\scriptstyle\oplus$}}{!}{$\boxplus$}}}}
      {\vcenter{\hbox{\resizebox{\widthof{$\scriptscriptstyle\oplus$}}{!}{$\boxplus$}}}}
  }\displaylimits
}

\begin{document}


\title[Distributed Persistent Homology for 2D Alpha Complexes]{Distributed Persistent Homology for 2D Alpha Complexes}


\author[F. Jensen]{Freya Jensen}\email{freya.jensen@iwr.uni-heidelberg.de}

\author[\'A. Torras-Casas]{\'Alvaro Torras-Casas}\email{atorras@us.es}

\address[F. Jensen]{IWR, Universit\"at Heidelberg, Im Neuenheimer Feld 225, Heidelberg, 69120, Germany}

\address[\'A. Torras-Casas]{Dpto. Matematica Aplicada I, Universidad de Sevilla, ETSII, Avenida Reina Mercedes s/n, Sevilla, 41012,  Spain}


\begin{abstract}
    We introduce a new algorithm to parallelise the computation of persistent homology of 2D alpha complexes. Our algorithm distributes the input point cloud among the cores which then compute a cover based on a rectilinear grid. 
    We show how to compute the persistence Mayer--Vietoris spectral sequence from these covers and how to obtain persistent homology from it. 
    For this, we introduce second-page collapse conditions and explain how to solve the extension problem. Finally, we give an overview of an implementation in C++ using Open MPI and discuss some experimental results.
\end{abstract}

\keywords{Persistent Homology; Distributed Computation; Alpha Complex; Mayer--Vietoris; Spectral Sequences}

\subjclass[Mathematics Subject Classification]{55N31, 55T99, 55-04}

\maketitle


\section{Introduction}\label{sec1}

We present an algorithm to compute persistent homology of a two-dimensional point cloud in parallel.
Our work combines an algorithm by Lo~\cite{Lo2012}, which constructs and stores the Delaunay triangulation of a two-dimensional point cloud in a distributed way, with an algorithm from~\cite{Torras2023}, which computes the persistent homology of a simplicial complex in parallel by using a finite cover. 
All in all, we pursue two goals with our approach: 
We want to compute the exact persistent homology of large data sets while gaining valuable additional information on the local persistent homology of the covering subcomplexes.

Before we describe our algorithm and its implementation in more detail, we give a short overview of related work.
First, the algorithm for the Delaunay triangulation described in~\cite{Lo2012} is a parallel incremental insertion algorithm based on a rectilinear grid.
There are other methods for distributed Delaunay triangulation construction in two dimensions.
An overview of several parallel algorithms may be found in~\cite{Kohout2005}.
One of the advantages of using the rectilinear grid-based algorithm by Lo is that it allows us to define the intersections between the covers in a straightforward way. 
Even though we cannot guarantee that the individual subcomplexes are simply connected, they are connected.
Then we use the Delaunay triangulation to compute the alpha complex.
That is, we compute a filtration value for each simplex in the complex.
An algorithm for the parallel computation of alpha complexes for biomolecules was presented in~\cite{Masood2020}.
This method does not require prior computation of the Delaunay triangulation
but it also does not provide us with a cover of the resulting complex.

Spectral sequences have been used for distributing computations of cohomology groups in~\cite{CuGhNa2016} and recently,
in~\cite{YoonGhrist2020}, spectral sequences are used for distributing persistent homology computations.
However, all of~\cite{CuGhNa2016, YoonGhrist2020} assume that the nerve of the cover is one-dimensional.
The algorithm presented by Lewis and Morozov in~\cite{LewisMorozov2015} is very similar to ours in its basic ideas. 
They use a spatial decomposition of the domain and build the Mayer-Vietoris blow-up complex of the underlying simplicial complex corresponding to the spatial decomposition.
Another similarity of~\cite{LewisMorozov2015} with our work is that we also use a processor for each cover element.
The fundamental difference to our approach is that Lewis and Morozov directly compute the persistent homology of the blow-up complex by parallelizing the matrix reduction.
In addition to~\cite{LewisMorozov2015}, there are algorithms which reduce the differential matrix for the original complex in parallel, such as~\cite{BauKerRei14, MoroNigm20}. In contraposition, our approach aims at distributing computations on both the filtered complex construction as well as persistent homology computation from a partition of the input dataset.

Besides the parallel algorithms, there are many sequential methods for computing the persistent homology of a point cloud. 
One of the fastest sequential algorithms is Ripser~\cite{Bauer2021Ripser}, which computes the Vietoris-Rips complex for the data and computes the corresponding persistent homology.
Another library for computing the alpha complex of a point cloud and its persistent homology is GUDHI~\cite{gudhi}. 

Our algorithm develops the method from~\cite{Lo2012} further so that each processor which handles one cover also computes all intersections with other covers. 
In addition, each processor also computes the Alpha filtration for the complexes it is handling.
Both these computations depend on additional information from the other covers and we aimed at minimising the necessary data exchange between the processors while keeping the computations exact.
Next, all processors compute persistent homology and merge their results by computing the second page of the Mayer-Vietoris spectral sequence as well as solving the extension problem; these last two steps are coordinated by two and one processors,respectively.

For all randomly generated data sets that we tried, we have observed that the spectral sequence collapses at the second page. 
This observation is important as it reduces the complexity of our algorithm substantially.
Unfortunately, one cannot generally assume that the spectral sequence collapses on the second page, as we show in Example~\ref{ex:second-page-differential}.
This is why we give conditions in Subsection~\ref{sub:second_collapse} as to when this is the case.

We have implemented our algorithm in C++ on top of MPI to show that the method is capable of dealing with large data sets consisting of around a million data points.
In addition, we use the \texttt{Simplex\_tree} 
structure implemented in GUDHI~\cite{gudhi} to store the simplicial complexes and the persistent homology computation implemented in PHAT~\cite{phat}, which is a library of several sequential algorithms for computing the persistent homology given the boundary matrix of a filtered complex.
Also, in our experiments, we check that we obtain the same persistent homology barcode (up to a tolerance) as using standard methods.

Now we describe the structure of the article. First, in Section~\ref{sec:background}, we give a short introduction to the Delaunay triangulaton $\DX$ and Alpha complexes. In this section we also introduce persistent homology, barcode bases and associated matrices to persistent morphisms. Next, in Section~\ref{sec:alpha-construction}, we start describing the construction of a covering set for $\DX$ based on Lo~\cite{Lo2012}. 
We move on to compute double and triple intersections of such covers, as well as the computation of their associated filtrations. In Section~\ref{sec:permaviss_alpha} we introduce the Persistence Mayer-Vietoris spectral sequence associated to the cover constructed in Section~\ref{sec:alpha-construction}. We explain how to compute such a spectral sequence and how to obtain persistent homology from its computation. 
Also, this description includes a discussion of when this spectral sequence collapses on the second page.
Afterwards, in Section~\ref{sec:implementation}, we give an overview of our implementation as well as some experiments. 
Finally, we conclude the article and discuss future research directions.

\section{Background}\label{sec:background}

\subsection{Delaunay Triangulations}
Before we introduce Delaunay triangulations, we provide the reader with basic notions on simplical complexes.

An (abstract) \textit{simplicial complex} $K$ on a set $X$ is a collection of subsets of $X$ (a.k.a simplices) such that for any element $\sigma \in K$ if $\tau \subseteq \sigma$ then $\tau \in K$. A simplex $\sigma \in K$ containing $n+1$ elements is called an \emph{$n$-simplex} and a subset $\tau \subseteq \sigma$ is called a \emph{face} of $\sigma$. 
Given a simplicial complex $K$ and a subset $L\subseteq K$, if $L$ is a simplicial complex then we call it a \emph{simplicial subcomplex (of $K$)} or, for short, a \emph{subcomplex (of $K$)}. 
A simplicial complex $K$ is said to be of dimension $n$ if all its simplices have at most $n+1$ elements. In addition, $0$-simplices, $1$-simplices and $2$-simplices are also called \emph{vertices}, \emph{edges} and \emph{triangles}, respectively.
We denote by $K^p$ the subset of $K$ composed of all $p$-simplices.
A \textit{filtered simplicial complex} is a finite collection of simplicial complexes $\{K_i\}_{i\in I}$, $I\subset \bR$ finite, such that whenever $i<j$, $K_i$ is a subcomplex of $K_j$. We often say that $\{K_i\}_{i\in I}$ \emph{defines a filtration} on $\bigcup_{i \in I} K_i$.

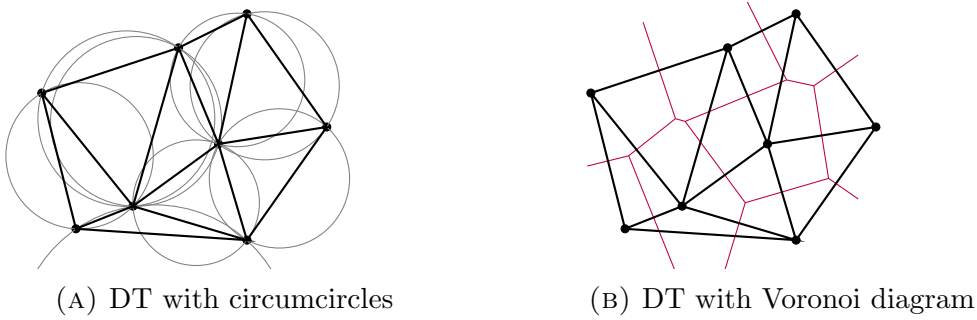
\begin{figure}
    \centering
    \begin{subfigure}[b]{0.45\textwidth}
    \centering
    \begin{tikzpicture}[scale = .75]
        \clip (-4.3,-1.5) rectangle (3.5,3.2);
        \coordinate (A) at (0,3);
        \coordinate (B) at (0,-1);
        \coordinate (C) at (1.4,1);
        \coordinate (D) at (-3.6,1.6);
        \coordinate (E) at (-2, -0.4);
        \coordinate (F) at (-0.5, 0.7);
        \coordinate (G) at (-1.2,2.4);
        \coordinate (H) at (-3,-0.8);
        \foreach \i in {A,B,C,D,E,F,G,H} {
            \filldraw (\i) circle (2pt);
        }
        \path let \p1=($(H)-(B)$),\p2=($(H)!0.5!(B)$),
            \p3=($(H)-(E)$),\p4=($(H)!0.5!(E)$), 
            \n1={(-(\x2*\x3)+\x3*\x4+\y3*(-\y2 +\y4))/(\x3*\y1-\x1*\y3)} in coordinate (centerHBE) at ({\x2+\n1*\y1},{\y2-\n1*\x1});
        \draw[circle through 3 points={A}{F}{G}, color=gray];
        \draw[circle through 3 points={F}{C}{A}, color=gray];
        \draw[circle through 3 points={F}{C}{B}, color=gray];
        \draw[circle through 3 points={F}{E}{B}, color=gray];
        \draw[circle through 3 points={F}{E}{G}, color=gray];
        \draw[circle through 3 points={D}{E}{G}, color=gray];
        \draw[circle through 3 points={H}{D}{E}, color=gray];
        \draw[half circle through 3 points={B}{H}{E}, color=gray];
        \draw[thick] (A) -- (F) -- (G) -- (A);
        \draw[thick] (A) -- (C) -- (F) -- (B) -- (C);
        \draw[thick] (E) -- (H) -- (D) -- (E) -- (G) -- (D);
        \draw[thick] (F) -- (E) -- (B) -- (H);
    \end{tikzpicture}
    \caption{DT with circumcircles}
    \label{fig:delaunay_triangulation_circle}
    \end{subfigure}
    \begin{subfigure}[b]{0.45\textwidth}
    \centering
     \begin{tikzpicture}[scale = .75]
        \clip (-4,-1.5) rectangle (3.5,3.2);
        \coordinate (A) at (0,3);
        \coordinate (B) at (0,-1);
        \coordinate (C) at (1.4,1);
        \coordinate (D) at (-3.6,1.6);
        \coordinate (E) at (-2, -0.4);
        \coordinate (F) at (-0.5, 0.7);
        \coordinate (G) at (-1.2,2.4);
        \coordinate (H) at (-3,-0.8);
        \foreach \i in {A,B,C,D,E,F,G,H} {
            \filldraw (\i) circle (2pt);
        }
        \path let \p1=($(A)-(F)$),\p2=($(A)!0.5!(F)$),
            \p3=($(A)-(G)$),\p4=($(A)!0.5!(G)$), 
            \n1={(-(\x2*\x3)+\x3*\x4+\y3*(-\y2 +\y4))/(\x3*\y1-\x1*\y3)} in coordinate (centerAFG) at ({\x2+\n1*\y1},{\y2-\n1*\x1});
        \path let \p1=($(A)-(F)$),\p2=($(A)!0.5!(F)$),
            \p3=($(A)-(C)$),\p4=($(A)!0.5!(C)$), 
            \n1={(-(\x2*\x3)+\x3*\x4+\y3*(-\y2 +\y4))/(\x3*\y1-\x1*\y3)} in coordinate (centerAFC) at ({\x2+\n1*\y1},{\y2-\n1*\x1});
        \path let \p1=($(B)-(F)$),\p2=($(B)!0.5!(F)$),
            \p3=($(B)-(C)$),\p4=($(B)!0.5!(C)$), 
            \n1={(-(\x2*\x3)+\x3*\x4+\y3*(-\y2 +\y4))/(\x3*\y1-\x1*\y3)} in coordinate (centerBFC) at ({\x2+\n1*\y1},{\y2-\n1*\x1});
        \path let \p1=($(B)-(F)$),\p2=($(B)!0.5!(F)$),
            \p3=($(B)-(E)$),\p4=($(B)!0.5!(E)$), 
            \n1={(-(\x2*\x3)+\x3*\x4+\y3*(-\y2 +\y4))/(\x3*\y1-\x1*\y3)} in coordinate (centerBFE) at ({\x2+\n1*\y1},{\y2-\n1*\x1});
        \path let \p1=($(G)-(F)$),\p2=($(G)!0.5!(F)$),
            \p3=($(G)-(E)$),\p4=($(G)!0.5!(E)$), 
            \n1={(-(\x2*\x3)+\x3*\x4+\y3*(-\y2 +\y4))/(\x3*\y1-\x1*\y3)} in coordinate (centerFEG) at ({\x2+\n1*\y1},{\y2-\n1*\x1});
        \path let \p1=($(G)-(D)$),\p2=($(G)!0.5!(D)$),
            \p3=($(G)-(E)$),\p4=($(G)!0.5!(E)$), 
            \n1={(-(\x2*\x3)+\x3*\x4+\y3*(-\y2 +\y4))/(\x3*\y1-\x1*\y3)} in coordinate (centerDEG) at ({\x2+\n1*\y1},{\y2-\n1*\x1});
        \path let \p1=($(H)-(D)$),\p2=($(H)!0.5!(D)$),
            \p3=($(H)-(E)$),\p4=($(H)!0.5!(E)$), 
            \n1={(-(\x2*\x3)+\x3*\x4+\y3*(-\y2 +\y4))/(\x3*\y1-\x1*\y3)} in coordinate (centerHDE) at ({\x2+\n1*\y1},{\y2-\n1*\x1});
        \path let \p1=($(H)-(B)$),\p2=($(H)!0.5!(B)$),
            \p3=($(H)-(E)$),\p4=($(H)!0.5!(E)$), 
            \n1={(-(\x2*\x3)+\x3*\x4+\y3*(-\y2 +\y4))/(\x3*\y1-\x1*\y3)} in coordinate (centerHBE) at ({\x2+\n1*\y1},{\y2-\n1*\x1});
        \draw[color=purple] (centerBFE) -- (centerBFC) -- (centerAFC) -- (centerAFG) -- (centerFEG) -- (centerDEG) -- (centerHDE) -- (centerHBE);
        \draw[color=purple] (centerHBE) -- (centerBFE) -- (centerFEG);
        \coordinate (AC) at ($(A)!0.5!(C)$);
        \coordinate (BC) at ($(B)!0.5!(C)$);
        \coordinate (AG) at ($(A)!0.5!(G)$);
        \coordinate (DG) at ($(D)!0.5!(G)$);
        \coordinate (DH) at ($(D)!0.5!(H)$);
        \draw[ray through 2 points={centerAFC}{AC}, color=purple];
        \draw[ray through 2 points={centerAFG}{AG}, color=purple];
        \draw[ray through 2 points={centerDEG}{DG}, color=purple];
        \draw[ray through 2 points={centerHDE}{DH}, color=purple];
        \draw[long ray through 2 points={centerBFC}{BC}, color=purple];
        \draw[thick] (A) -- (F) -- (G) -- (A);
        \draw[thick] (A) -- (C) -- (F) -- (B) -- (C);
        \draw[thick] (E) -- (H) -- (D) -- (E) -- (G) -- (D);
        \draw[thick] (F) -- (E) -- (B) -- (H);
    \end{tikzpicture}
    \caption{DT with Voronoi diagram}
    \label{fig:delaunay_triangulation_voronoi}
    \end{subfigure}
    \caption{Delaunay triangulation (DT) of some random points}
    \label{fig:delaunay_triangulation}
\end{figure}

Let $\bX$ be a finite set of points in $\bR^2$.
We consider a two-dimensional simplicial complex on $\bX$ such that its vertices are the points from $\bX$, its edges span segments between two different points from $\bX$ and its triangles are composed of three non-colinear points; this is also known as a \emph{geometric simplicial complex} on $\bX$. In particular,
a \textit{triangulation} of $\bX$ is a two-dimensional geometric simplicial complex $\TX$ such that
$\bX$ is the set of vertices in $\TX$, and the union of all simplices in $\TX$ is the convex hull of $\bX$.
A triangle on $\bX$ is called \textit{Delaunay} with respect to $\bX$ if no points of $\bX$ are lying in its open circumcircle.
If all the triangles from $\TX$ are Delaunay then $\TX$ is called a \emph{Delaunay Triangulation} (of $\bX$).
Such a Delaunay triangulation is depicted in Figure~\ref{fig:delaunay_triangulation_circle}, where the circumcircles of its triangles are shown in gray.
One can show that a Delaunay triangulation $\DX$ for $\bX$ always exists; however, in general, it is not unique.

\leavevmode

Apart from the vertices of a triangle, other points of $\bX$ are allowed to lie on the boundary of the corresponding circumcircle as well.
In this case, there are three ways to define a Delaunay structure on the underlying set $\bX$.
The first option is to include all Delaunay triangles.
Then the resulting simplicial complex may not be embedded into $\bR^2$ anymore.
The second option is to choose a minimal subset of Delaunay triangles such that the resulting simplicial complex is a proper triangulation of $\bX$.
This choice is not unique.
The last option is to exclude all crossing edges inside the circumcircle of the cocircular vertices.
This option does not result in a proper triangulation.
If the points of $\bX$ are in \textit{general position} though, i.e. no four points are lying on the same circle, the Delaunay triangulation $\DX$ is unique.

Delauanay himself provided a very helpful condition to check whether a triangulation $\TX$ is Delaunay or not:
\begin{lemma}[Lemma of Delaunay]\label{thm:delaunay}
    Let $\TX$ be a triangulation of the finite point set $\bX$. Then a triangle of $\TX$ is Delaunay if and only if none of the vertices of its adjacent triangles lies inside its circumcircle.
\end{lemma}

Alternatively, one can also define the Delaunay triangulation of $\bX$ as the dual of the Voronoi diagram of $\bX$.
The \textit{Voronoi cell} $V_x$ of a point $x$ in $\bX$ is the set of all points $y\in\bR^2$ that are closer to $x$ or equally close to $x$ as to any other point of $\bX$.
The collection $\VorX = \{V_x \vert x \in \bX\}$ of all Voronoi cells is called \textit{Voronoi diagram} of $\bX$.

Let $I$ be an index set. Given a cover $\cU=\{U_{\alpha}\}_{\alpha\in I}$ of $\bR^2$,
we define a simplicial complex
$\cN(\cU)$---called the nerve of $\cU$---as follows:
$\sigma = (\alpha_0,\ldots,\alpha_n)$ is an $n$-simplex of $\cN(\cU)$ if and only if the intersection of the sets $U_{\alpha_0},\ldots,U_{\alpha_n}$ is non empty, i.e.
\[
\cN(\cU) = \left\{\sigma \subset I \Bigg\vert \bigcap_{\alpha \in \sigma} U_{\alpha} \neq \emptyset\right\}.
\]

As an abstract simplicial complex, the Delaunay triangulation $\DX$ is isomorphic to the nerve complex of the Voronoi diagram $\VorX$ of $\bX$.
The Voronoi cells of the Delaunay triangulation in Figure~\ref{fig:delaunay_triangulation} are highlighted in red in Figure~\ref{fig:delaunay_triangulation_voronoi}. Later we use that $\DX$ is contractible since it is collapsible~\cite{AttaliLieutierSalinas2019}.

\subsection{Alpha Complex}\label{sec:alpha_complex}
Starting with the Voronoi diagram $\VorX$, we can create a family of spaces together with corresponding covers:
For any non-negative real number $r$, consider the union of closed balls of diameter $r$ around each point $x\in\bX$, $\bigcup_{x\in\bX} B_x(r)$.
We can construct a cover for this union based on intersections between 
closed balls with the corresponding Voronoi cells, $R_x(r) = B_x(r) \cap V_x$.
We define the \textit{Alpha complex} as the nerve of this cover, i.e.
\[\AX{r} = \left\{ \sigma \subset \bX \Bigg\vert \bigcap_{x\in\sigma} R_x(r) \neq \emptyset\right\}.\]
For any $r>0$, $\AX{r}$ is a subcomplex of $\DX$ as $R_x(r)\subset V_x$ for any $x\in\bX$.
The above approach leads to an important observation. Let $\vert K\vert$ be the \emph{geometric realisation} of a simplicial complex $K$.

\begin{theorem}[Nerve theorem]
    Let $\cU$ be a finite collection of closed, convex sets in Euclidean space.
    Then, $\vert \cN(\cU)\vert$ and the union of $\cU$ are homotopy equivalent.
\end{theorem}

Since $\{R_x(r)\}_{x\in\bX}$ is a cover of $\bigcup_{x\in\bX} B_x(r)$ and all $R_x(r)$ are closed and convex, by the Nerve theorem we get a homotopy equivalence
\[ \vert \AX{r} \vert \simeq \bigcup_{x\in\bX} B_x(r).\]

Since, for any $r<s$, $\AX{r}$ is a subcomplex of $\AX{s}$ and there is an $r^\ast$ such that $\AX{r} = \DX$ for all $r\geq r^\ast$,
the collection $\{\AX{r}\}_{r>0}$ defines a filtration on $\DX$.

Alternatively, we can construct the Alpha complex by defining a filtration on the Delaunay triangulation.
In order to do this, we differentiate between two kinds of edges, Gabriel and non-Gabriel.
An edge is \textit{Gabriel} if the interior of its circumcircle is empty.
We define a filtration function, f: $\DX \rightarrow \bR$, by assigning to each triangle as well as to each Gabriel edge the radius of its circumcircle.
In the case of a non-Gabriel edge, we assign to it the minimum filtration value of its adjacent triangles.
To any vertex we assign the value $0$.
One can easily see that the complexes defined by the resulting sublevel sets are equal to the corresponding Alpha complexes, i.e. $f^{-1}(-\infty,r) = \AX{r}$.

\begin{example}
Consider the Delaunay triangulation on four points illustrated in Figure~\ref{fig:non_gabriel_edge}, where the vertical edge is non-Gabriel.
In this case, the Voronoi cells corresponding to the vertices of the edge do not intersect at the smallest distance between the two points, i.e. at the radius of its circumcircle.
This intersection point is depicted in green in the figure.
The Voronoi cell of the point making the edge non-Gabriel, 
the point on the right in Figure~\ref{fig:non_gabriel_edge}, is causing the two other Voronoi cells to intersect at a greater distance.
Namely, at the radius of the circumcircle of the triangle
having the three of them as vertices.
The center of this circumcircle is depicted as a red dot in Figure~\ref{fig:non_gabriel_edge}.
The boundaries of the Voronoi cells are depicted in red as well.
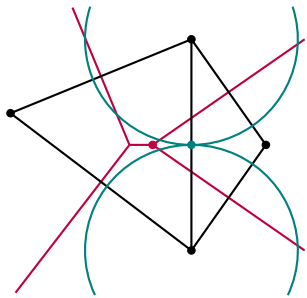
\begin{figure}[ht]
    \centering
    \begin{tikzpicture}[scale = .7]
        \coordinate (A) at (0,3);
        \coordinate (B) at (0,-1);
        \coordinate (C) at (1.4,1);
        \coordinate (D) at (-3.4,1.6);
        \coordinate (E) at ($(A)!0.5!(B)$);
        \path let \p1=($(A)-(B)$),\p2=($(A)!0.5!(B)$),
            \p3=($(A)-(C)$),\p4=($(A)!0.5!(C)$), 
            \n1={(-(\x2*\x3)+\x3*\x4+\y3*(-\y2 +\y4))/(\x3*\y1-\x1*\y3)} in coordinate (F) at ({\x2+\n1*\y1},{\y2-\n1*\x1});
        \path let \p1=($(A)-(B)$),\p2=($(A)!0.5!(B)$),
            \p3=($(A)-(D)$),\p4=($(A)!0.5!(D)$), 
            \n1={(-(\x2*\x3)+\x3*\x4+\y3*(-\y2 +\y4))/(\x3*\y1-\x1*\y3)} in coordinate (G) at ({\x2+\n1*\y1},{\y2-\n1*\x1});
        \foreach \i in {A,B,C,D} {
            \filldraw (\i) circle (2pt);
        }
        \coordinate (auxAC) at ($(A)!0.5!(C)$);
        \coordinate (auxBC) at ($(B)!0.5!(C)$);
        \draw[ray through 2 points={F}{auxAC}, color=purple, thick];
        \draw[ray through 2 points={F}{auxBC}, color=purple, thick];
        \coordinate (auxAD) at ($(A)!0.5!(D)$);
        \coordinate (auxBD) at ($(B)!0.5!(D)$);
        \draw[ray through 2 points={G}{auxAD}, color=purple, thick];
        \draw[long ray through 2 points={G}{auxBD}, color=purple, thick];
        \draw[thick, color=purple] (F) -- (G);
        \draw[thick] (A) -- (C) -- (B) -- (A) -- (D) -- (B);
        \filldraw[color=teal] (E) circle (2pt);
        \draw[color=teal, thick] (A) ++(18:2.0) arc(18:-198:2.0);
        \draw[color=teal, thick] (B) ++(-25:2.0) arc(-25:205:2.0);
        \filldraw[color=purple] (F) circle (2pt);
    \end{tikzpicture}
    \caption{Shortest distance (green) between two vertices of a non-Gabriel edge versus actual point of first intersection (red) of the corresponding Voronoi cells}
    \label{fig:non_gabriel_edge}
\end{figure}
\end{example}

\begin{figure}
    \centering
    \includegraphics[width=0.9\textwidth]{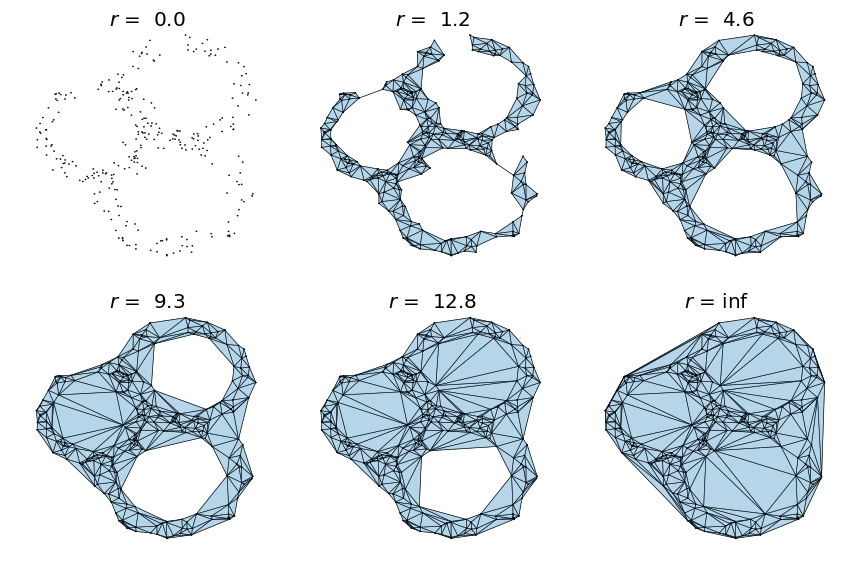}
    \caption{Depiction of $\AX{r}$ for increasing values.}
    \label{fig:alpha_filtered}
\end{figure}

\begin{example}\label{ex:alpha_filtered}
Consider the point cloud $\bX$ depicted on the top left of Figure~\ref{fig:alpha_filtered}. Within the same figure, we depict an increasing sequence of alpha complexes $\AX{r}$ where $r$ ranges over the set $R=\{0, 1.2, 4.6, 9.3, 12.8, \infty\}$.
\end{example}
\leavevmode

\subsection{Persistent Homology of Alpha Complexes}\label{sub:PH-alpha}
We want to investigate the global structure of the point set $\bX$.
To do this, we use the field $\bZ_2$ and we use \emph{simplicial homology}; for an introduction to this concept, see~\cite{roadmap}. Our choice of field is because we work with $2$ dimensional alpha complexes embedded in the plane $\bR^2$, and so without torsion. Fixing $r \geq 0$, we consider $\AX{r}$ as an abstract simplicial complex
and represent its simplices as tuples on $\bX$ that are equivalent under permutations. Thus, fixed an integer $m\geq 0$, we write $C_m(\AX{r})=\bZ_2[\AX{r}]$ which is the vector space spanned by elements from the aforementioned representation of $\AX{r}$. 
The elements from $C_m(\AX{r})$ are called \emph{$m$-chains} (on $\AX{r}$). We also denote as  $d_m:C_m(\AX{r})\rightarrow C_{m-1}(\AX{r})$ the \emph{$m$-differential}.
This leads us to the simplicial homology groups $\Ho_m(\AX{r}; \bZ_2)= \Ker(d_m)\big/ \Ima(d_{m+1})$. For short, we write $\Ho_m(\AX{r})$ instead of $\Ho_m(\AX{r}; \bZ_2)$.

To introduce persistent homology, we consider variations over the parameter $r\geq 0$.
Since $\bX$ is finite, there is a finite sequence 
\[\emptyset \subseteq \AX{r_0} \subseteq \AX{r_1} \subseteq \ldots \subseteq \AX{r_N} = \DX,\]
such that $\AX{r}$ is constant along $r \in \left[r_n,r_{n+1}\right)$ for $0 \leq n \leq N-1$.
For simplicity, we write $\AX{n}$ for $\AX{r_n}$ in the following. 
Then, the inclusions induce a sequence of finite-dimensional $\bZ_2$-vector spaces
\[0 \rightarrow \Ho_m(\AX{0}) \rightarrow \Ho_m(\AX{1}) \rightarrow \ldots \rightarrow \Ho_m(\AX{N}) = \Ho_m(\DX).\]
The above sequence is the \textit{persistent homology} of $\AXn$ at dimension $m\geq 0$, and is denoted as $\PH_m(\AXn)$. We also call $\PH_m(\AXn)$ the \textit{alpha complex persistence} of $\bX$.
\begin{example}\label{ex:PH_alpha}
    Consider the alpha complexes $\AX{r_i}$ with $r_i \in R$ from Example~\ref{ex:alpha_filtered}. Now, we compute the homology of these complexes (in dimension $1$) as well as the structure maps for $r_i \in R$. This leads to the sequence
    \begin{equation}\label{seq:ex_barcode_alpha}
    \begin{tikzcd}[/tikz/column 6/.style={column sep=-0.5em}]
    0 \ar[r] & 
    \bZ_2 \oplus \bZ_2 
    \ar[r, "\sbm{
        1 \, 0 \\
        0 \, 1 \\
        0 \, 0 
    }"] 
    &
    \bZ_2 \oplus \bZ_2  \oplus \bZ_2 
    \ar[r, "\sbm{
        1 \, 0 \, 0\\
        0 \, 0 \, 1
    }"]
    &
    \bZ_2 \oplus \bZ_2 
    \ar[r, "\sbm{
        1 \, 0 
    }"]
    &
    \bZ_2 \ar[r]
    & 
    0
    &.
    \end{tikzcd}
    \end{equation}
    Notice that Sequence~(\ref{seq:ex_barcode_alpha}) tracks the nontrivial cycles depicted in Figure~\ref{fig:alpha_filtered}. For example, the second term from~(\ref{seq:ex_barcode_alpha}), which corresponds to $r_1=1.2$, has two copies of $\bZ_2$ that indicate the two nontrivial cycles on the middle top complex from Figure~\ref{fig:alpha_filtered}.
\end{example}
\begin{definition}
A \emph{persistence module} $\bV$ (over $\bZ_2$) is a collection of $\bZ_2$-vector spaces $\bV_t$, for all $t \in \bR$, together with linear maps (a.k.a. \emph{structure maps}) $\rho_{st}:\bV_s\rightarrow \bV_t$ for all pairs $s \leq t$ such that: 1) $\rho_{tt}$ is the identity on $\bV_t$ and 2) given $r \leq s \leq t$, the equality $\rho_{st}\rho_{rs}=\rho_{rt}$ holds.
\end{definition}
We also refer to a persistence module as a pair $(\bV, \rho)$ to indicate that $\rho$ denotes the structure maps of $\bV$. 
Given two persistence modules $(\bV, \rho)$ and $(\bW, \tau)$, a morphism $f$ between $\bV$ and $\bW$ (a.k.a. a \emph{persistence morphism}), denoted as $f:\bV\rightarrow \bW$, is a collection of linear maps $f_t:\bV_t\rightarrow \bW_t$ such that, for any pair $s\leq t$, the equality $\tau_{st}f_s = f_t\rho_{st}$ holds. 

\begin{figure}
    \centering
    \includegraphics[width=0.8\textwidth, height=0.25\textwidth]{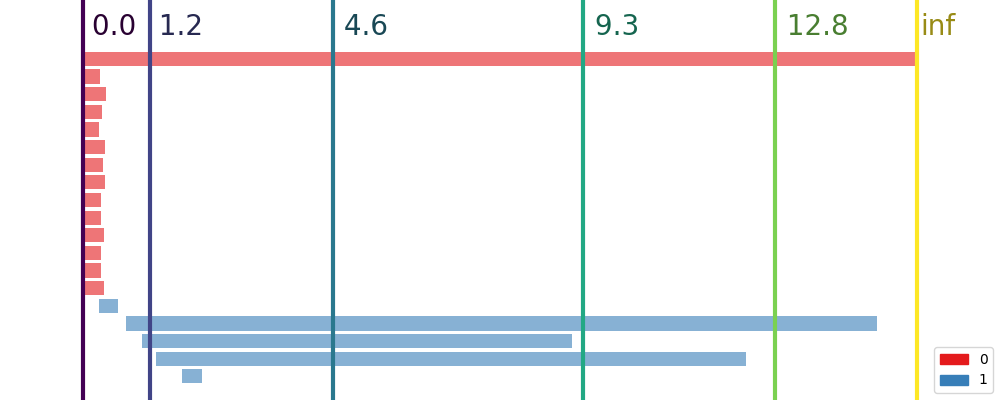}
    \caption{Persistent homology barcodes in dimension 0 (red) and in dimension 1 (blue) of $\AXn$. One red bar ``never ends'', which we indicate with a vertical yellow line with the ``inf'' label. }
    \label{fig:PH_alpha_vert}
\end{figure}

Under very general conditions~\cite{CrawleyBoevey2015}, a persistence module can be decomposed into very simple pieces. Given $a\leq b$, an \emph{interval module} $\bZ_{2[a,b)}$ is a persistence module such that $\bZ_{2[a,b)t}=\bZ_2$ for all $t \in [a,b)$ and its structure maps are either the identity whenever possible or zero otherwise. A persistence module $\bV$ is \emph{tame} if there exists an isomorphism $\bV\simeq \bigoplus_{i=1}^N \bZ_{2[a_i, b_i)}$; i.e. $\bV$ is isomorphic to the direct sum of finitely many interval modules. In this article, all persistence modules are tame.
The multiset of intervals in the decomposition of $\bV$ is known as the \emph{barcode of $\bV$} and is denoted by $\rB(\bV)$.
\begin{example}~\label{ex:barcode_alpha}
    Computing directly Sequence~(\ref{seq:ex_barcode_alpha}) from Example~\ref{ex:PH_alpha} can be a tedious task. This is why the standard procedure computes the barcode of $\PH_m(\AXn)$, which encodes all data in $\PH_m(\AXn)$ very compactly. Figure~\ref{fig:PH_alpha_vert} depicts $\rB(\PH_m(\AXn))$ for dimensions $m=0,1$, as well as the values from $R$ indicated by vertical lines.
\end{example}
\leavevmode

\subsection{Barcode bases}
In this subsection, we briefly review barcode bases and refer to ~\cite{Torras2023} for a more detailed exposition; other articles such as~\cite{JacNanTil23, YooGhrGiu2023} also study barcode bases.
\begin{definition}\label{def:barcode-basis}
    Let $\bV$ be a tame persistence module. A \emph{barcode basis} $\cA$ of $\bV$ is an isomorphism $\alpha: \bigoplus_{i=1}^N \bZ_{2[a_i, b_i)}\rightarrow \bV$. A \emph{generator} is the restriction of $\alpha$ to an interval module $\alpha_i:\bZ_{2[a_i, b_i)}\rightarrow \bV$. We also refer to a barcode basis as a set $\cA=\{\alpha_i\}_{i=1}^N$.
\end{definition}

\begin{figure}
    \centering
    \begin{subfigure}[b]{\textwidth}
        \centering
        \includegraphics[width=\textwidth]{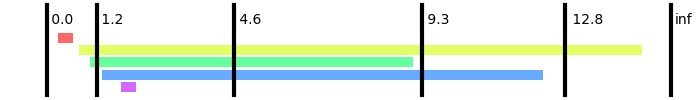}
        \caption{Barcode of $\PH_1(\AXn)$ together with values from $R$ as vertical lines.}
        \label{subfig:PH_1D_alpha}
    \end{subfigure}
    \vfill 
    \begin{subfigure}[b]{\textwidth}
        \centering
        \includegraphics[width=\textwidth]{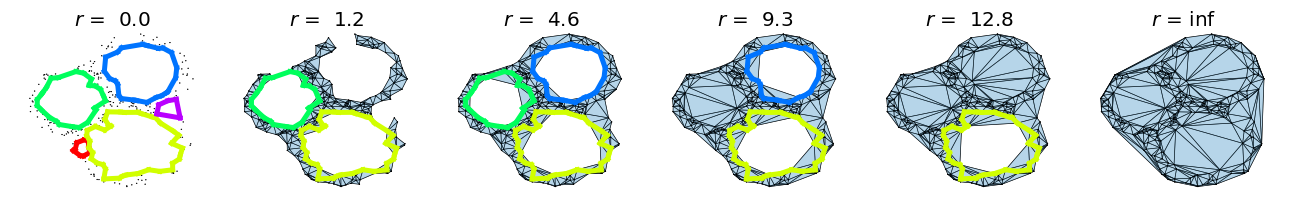}
        \caption{The complex $\AX{r_i}$ for $r_i$ ranging over $R$ together with nontrivial cycle representatives. We depict all representatives at filtration $0$ for illustration purposes. }
        \label{subfig:reps_1D_alpha}
    \end{subfigure}
    \caption{Barcode of $\PH_1(\AXn)$ and cycle representatives using matching colours.}
    \label{fig:PH_1_alpha_reps}
\end{figure}

Following Definition~\ref{def:barcode-basis}, given a generator $\alpha_i:\bZ_{2[a_i, b_i)}\rightarrow \bV$, we say that $\alpha_i$ is associated to the interval $[a_i, b_i)$ and denote this as $\alpha_i\sim [a_i, b_i)$; here $a_i$ and $b_i$ are the \emph{birth} and \emph{death} values of $\alpha_i$ respectively. 
Notice that $\cA$ is such that, given $t \in \bR$, the set 
$\big\{ \alpha_{it}(1_{\bZ_2}) \mid t \in [a_i,b_i)\big\}$
is a basis for $\bV_t$.
\begin{definition}\label{def:basis-orders}
    A barcode basis $\cA=\{\alpha_i \sim [a_i, b_i)\}_{i=1}^N$ is sorted following the \emph{standard} order if, for $i<j$, either $a_i < a_j$ or,  $a_i=a_j$ and $b_i \geq b_j$. Also, we say that $\cA$ follows the  \emph{endpoint} order if, for $i<j$, either $b_i < b_j$ or, $b_i=b_j$ and $a_i \leq a_j$.
\end{definition}
A barcode basis $\cA=\{\alpha_i\sim [a_i,b_i)\}_{i=1}^N$ is closely related to persistent homology \emph{representatives}. 
That is, when computing the barcode of $\PH_m(\AXn)$ one can obtain, for all $i=1,\ldots,N$, cycles $w_i \in C_m(\AX{a_i})$ (i.e. $d_m(w_i)=0$) such that they determine a barcode basis $\cA$ for $\PH_m(\AXn)$. 
More concretely, we have $\alpha_{ia_i}(1_{\bZ_2}) = w_i + \Ima(d_{m+1})_{a_i}$ for all $i=1,\ldots, N$. 
\begin{example}\label{ex:barcode_basis}
    Let $\AX{r_i}$ with $r_i \in R$ from Example~\ref{ex:alpha_filtered}.
    We depict the barcode of $\PH_1(\AXn)$ together with the corresponding cycle representatives in Figure~\ref{fig:PH_1_alpha_reps}; these determine a barcode basis $\cB=\{\beta_i\}_{i=1}^5$ which is sorted in standard order; following the intervals---from top to bottom---from Subfigure~\ref{subfig:PH_1D_alpha}. Given  $\beta_4$ (represented by the blue circle on top of Subfigure~\ref{subfig:reps_1D_alpha}) we might restrict $\beta_4$ to $R$ and obtain the commutative diagram
    \[
    \begin{tikzcd}[column sep=0.8cm, /tikz/column 6/.style={column sep=-0.5em}]
     0 \ar[d] \ar[r]
     & 
     0 \ar[d] \ar[r]
     &
     \bZ_2 \ar[d, 
        "\sbm{
            0 \\
            0 \\
            1
        }"]
    \ar[r]
     &
     \bZ_2 \ar[d, 
        "\sbm{
            0 \\
            1
        }"] 
    \ar[r]
     &
     0 \ar[d] \ar[r]
     &
     0 \ar[d]
    \\ 
    0 \ar[r] & 
    \bZ_2 \oplus \bZ_2 
    \ar[r, "\sbm{
        1 \, 0 \\
        0 \, 1 \\
        0 \, 0 
    }"] 
    &
    \bZ_2 \oplus \bZ_2  \oplus \bZ_2 
    \ar[r, "\sbm{
        1 \, 0 \, 0\\
        0 \, 0 \, 1
    }"]
    &
    \bZ_2 \oplus \bZ_2 
    \ar[r, "\sbm{
        1 \, 0 
    }"]
    &
    \bZ_2 \ar[r]
    & 
    0
    &,
    \end{tikzcd}
    \]
    where we have used that $[a_4, b_4)=[1.3, 12.3)$ and Sequence~(\ref{seq:ex_barcode_alpha}) from Example~\ref{ex:PH_alpha}.
\end{example}
\leavevmode

A persistence vector in $\bV$, $\gamma$, is a morphism from an interval module into $\bV$; i.e. $\gamma:\bZ_{2[a,b)}\rightarrow \bV$. We denote by $\PVect(\bV)$ the set of persistence vectors in $\bV$. The \emph{barcode sum} is an operation $\boxplus:\PVect(\bV)\times\PVect(\bV)\rightarrow \PVect(\bV)$ which is such that, given $\gamma, \tau \in \PVect(\bV)$ with $\gamma \sim [a_\gamma,b_\gamma)$ and $\tau \sim [a_\tau, b_\tau)$, then $(\gamma \boxplus \tau)_t = \gamma_t + \tau_t$ for all $t \in [a_\gamma,b_\gamma)\cap[a_\tau, b_\tau)$; the interested reader can consult the details in~\cite[Def.1.14]{Torras2023}.

We end this section with some additional terminology.
Let $s \in \bR$. We define an operator $\bOne_s:\PVect(\bV)\rightarrow \PVect(\bV)$ such that, given $\gamma \in \PVect(\bV)$ with $\gamma\sim [a_\gamma,b_\gamma)$, $\bOne_s(\gamma)_r = \gamma_r$ for all $r \in [\max(s,a_\gamma), b_\gamma)$ and is zero otherwise.
Given a persistence morphism $f:\bV\rightarrow \bW$ and given $\gamma \in \PVect(\bV)$ with $\gamma\sim [a_\gamma,b_\gamma)$, we define $f(\gamma):\bZ_{2[a_\gamma, b_{f(\gamma)})} \rightarrow \bW$ where $b_{f(\gamma)}=\sup\big\{ r \in [a_\gamma, b_\gamma) \colon f_r \gamma_r \neq 0 \big\}$ and also $f(\gamma)_r = f_r\gamma_r$ for all $r \in [a_\gamma, b_{f(\gamma)})$. 

\subsection{Persistence morphisms and barcode bases}\label{sub:persistence-morphisms}
Consider a pair of barcode bases $\cA$ and $\cB$ for $\bV$ and $\bW$ respectively, and let $f:\bV\rightarrow \bW$ be a persistence morphism. As in linear algebra, there exists a matrix $F=(k_{\alpha,\beta})_{\alpha\in\cA, \beta\in\cB}$ associated to $f$. (In Appendix~\ref{app:associated-matrix} we explain how to obtain such matrix in practice.) This matrix is such that, given $\alpha \in \cA$ with $\alpha\sim [a_\alpha,b_\alpha)$, there exists a unique subset $\cS_\alpha \subseteq \cB$ such that 
$
f(\alpha)= \bOne_{a_\alpha}\Big(\bigBoxPlus_{\beta \in \cS_\alpha} k_{\alpha,\beta}\beta\Big)\, 
$; where $k_{\alpha,\beta}=0$ for all $\beta \in \cB\setminus \cS_\alpha$. It is worth noticing that, given $\alpha \in \cA$ with $\alpha\sim [a_\alpha, b_\alpha)$ and $\beta \in \cB$ with $\beta\sim [a_\beta, b_\beta)$, if $k_{\alpha,\beta}\neq 0$ then $a_\beta \leq a_\alpha \leq b_\beta \leq b_\alpha$. 

\begin{figure}
    \begin{minipage}{.2\textwidth}
        \centering
        \begin{subfigure}[b]{\textwidth}
        \includegraphics[width=\textwidth]{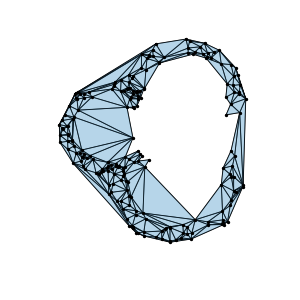}
            
            \caption{The subcomplex $K$.}
            \label{subfig:subcomplex}
        \end{subfigure}
    \end{minipage}
    \begin{minipage}{.7\textwidth}
        \centering
        \begin{subfigure}[b]{\textwidth}
        \includegraphics[width=\textwidth]{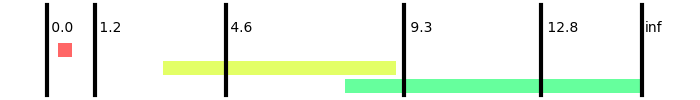}
            
            \caption{Barcode of $\PH_1(\AXn\cap K)$.}
            \label{subfig:barcode_subcomplex}
        \end{subfigure} 
        \begin{subfigure}[b]{\textwidth}
            \includegraphics[width=\textwidth]{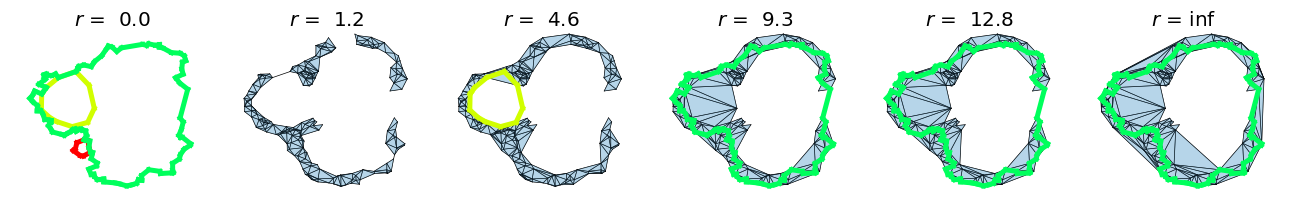}
            \caption{Representatives of persistent homology classes.}
            \label{subfig:reps_subcomplex}
        \end{subfigure} 
    \end{minipage}
    \caption{Subcomplex and its persistent homology barcode and representatives.}
\end{figure}

\begin{example}\label{ex:associated-matrix}
    Here we consider again Example~\ref{ex:barcode_basis} where $\cB=\{\beta_i\}_{i=1}^5$ is a barcode basis for $\PH_1(\AXn)$. We take a subcomplex $K$ from $\DX$ as shown in Subfigure~\ref{subfig:subcomplex}. This leads to the filtered complex $\AXn\cap K$ and to $\PH_1(\AXn\cap K)$; whose barcodes and representatives are depicted in Subfigures~\ref{subfig:barcode_subcomplex} and~\ref{subfig:reps_subcomplex} respectively. 
    We denote by $\cA=\{\alpha_i\}_{i=1}^3$ the barcode basis of $\PH_1(\AXn\cap K)$, where the generators are sorted following the bars from top to bottom in Subfigure~\ref{subfig:barcode_subcomplex}.
    One can then obtain 
    \[ 
    F=
    {\footnotesize
    \left[\begin{array}{ccc}
    1 & 0 & 0 \\
    0 & 0 & 1 \\
    0 & 1 & 1 \\
    0 & 0 & 1 \\
    0 & 0 & 0
    \end{array}\right]
    }
    \hspace{0.3cm}
    \mbox{ associated to }
    f:\PH_1(\AXn\cap K) \longrightarrow \PH_1(\AXn)\ ,
    \]
    in the bases $\cA$ and $\cB$. Intuitively, one can see that the representatives of $\alpha_1$ and $\alpha_2$ depicted as red and yellow cycles in Subfigure~\ref{subfig:reps_subcomplex} are equivalent (up to boundaries) to the representatives of $\beta_1$ and $\beta_3$ depicted as red and green cycles in Subfigure~\ref{subfig:reps_1D_alpha}. This is why the first two columns from $F$ are trivial except on the entries $F_{1,1}$ and $F_{3,1}$. Analogously, the representative of $\alpha_3$ is equivalent to the addition of the representatives of $\beta_2$, $\beta_3$, and $\beta_4$, as hinted by the third column from $F$.
\end{example}
\leavevmode

Next, we see why the barcode basis orders from Definition~\ref{def:basis-orders} are important in this work. Suppose that $\cA=\{\alpha_i\sim [a_i,b_i)\}_{i=1}^N$ is sorted following the standard order and $\cB=\{\beta_i\}_{i=1}^M$ is sorted following the endpoint order. Then, consider the associated matrix $F$ of $f$ and its Gaussian column reduction $R$. It turns out that one gets a barcode basis for $f(\bV)$ from the matrix $R$. More precisely, for each $\alpha_i \in \cA$, we consider the set $\cS_i=\{\beta_j \in \cB \colon R_{j,i}\neq 0\}$ and define the barcode basis of $f(\bV)$ given by
$
\cI = \Big\{ 
\bOne_{a_i}\Big( \bigBoxPlus_{\beta_j \in \cS_i} R_{j,i} \beta_j \Big)
\Big\}_{i=1}^N\,.
$
\begin{example}
    We continue from Example~\ref{ex:associated-matrix}. In this case, the columns from $F$ are already sorted in standard order, but the rows need to be sorted in endpoint order. From Subfigure~\ref{subfig:PH_1D_alpha}, we observe that in the endpoint order $\beta_1 < \beta_5 < \beta_3 < \beta_4 < \beta_2$. Thus, by reordering the rows from $F$, we obtain
    \[
    R = 
    {
    \footnotesize
    \left[\begin{array}{ccc}
    1 & 0 & 0 \\
    0 & 0 & 0 \\
    0 & 1 & 1 \\
    0 & 0 & 1 \\
    0 & 0 & 1 
    \end{array}\right]
    }
    \mbox{, which is already reduced. }
    \]
    Thus, denoting by $c_i$ the startpoints of the intervals associated to $\beta_i$ for all $i=1,2,3$, a barcode basis for $\Ima(f)$ is given by the set 
    $
    \{\bOne_{c_1}(\beta_1), 
    \bOne_{c_2}(\beta_3),
    \bOne_{c_3}(\beta_3 \boxplus \beta_4 \boxplus \beta_5)
    \}$. The intervals from $\rB(\Ima(f))$ are given by the birth of each column paired with the death of their pivots. In our case, the three intervals of the image barcode are $[c_1, b_1)$, $[c_2, b_3)$ and $[c_3, b_2)$, since the pivots from $R$ are $\beta_1$, $\beta_3$ and $\beta_2$ respectively.
\end{example}
One might track the reductions performed on $R$ to obtain a matrix of preimages $P$; which is a matrix indexed by $\cA$ in both the rows and columns. This matrix $P$ can be reordered by following the standard order on the columns while the endpoint order on the rows. Reducing this $P$ one obtains a barcode basis for $\Ker(f)$. More details can be found in~\cite{Torras2023}.

To work with spectral sequences, we often need to compute quotients of persistence modules. 
To compute quotients, we use a special type of Gaussian reduction introduced in~\cite{Torras2023}. In Appendix~\ref{app:example-quotient} we compute an explicit example of a quotient. 

\section{Alpha Filtration covers from Grids}\label{sec:alpha-construction}
In this section, we summarise how to obtain a (filtered) cover for the Delaunay triangulation based on a regular grid. 
We will then describe the computation of the resulting intersections in more detail.

\subsection{A cover for the Delaunay Triangulation}\label{sub:definition_cover}
In this subsection, we describe how the covers are defined and give some intuition for the resulting intersections.

Given a finite point-cloud $\bX$ embedded in $\bR^2$,
we follow the algorithm introduced in~\cite{Lo2012} to compute the Delaunay triangulation $\DX$ in parallel.
Since the point cloud is finite, 
it is contained within a compact, rectangular subspace $B$ of $\bR$ 
which we call the \emph{bounding box} of $\bX$.
We choose $B$ parallel to the $x$- and $y$-axis of the underlying coordinate system and such that none of the points in $\bX$ lie on the boundary of $B$.
Thus, we can describe $B$ by the minimum and maximum values in both $x$ and $y$ axes.
In order to compute the Delaunay triangulation in parallel,
$B$ is subdivided into $M$ rectangular zones $B_i$, for $i\in \{0,\ldots, M-1\}$, based on a rectilinear grid. 
Similarly to $B$, each of these smaller rectangular zones can again be described by its minimum and maximum values in both $x$ and $y$ axes,
$x^i_{\text{min}}, x^i_{\text{max}}, y^i_{\text{min}}, y^i_{\text{max}}$.
By using these zones, we can divide $\bX$ into disjoint subsets $\bX_i = \bX\cap B_i$, for $i\in \{0,\ldots, M-1\}$.
For points lying on the boundaries of zones, we make the following choice:
Let $p=(x_p,y_p)$ be such a point. 
Then $p$ belongs to $\bX_i$ if $x^i_{\text{min}} \leq x_p < x^i_{\text{max}}$
and if $y^i_{\text{min}} \leq y_p < y^i_{\text{max}}$.
We can then classify triangles of $\DX$ with respect to a fixed zone $B_i$.
A triangle from $\DX$ is called
\begin{itemize}
\item a \emph{boundary triangle of $B_i$} if some, but not all, of its vertices lie in $\bX_i$. 
\item an \emph{inner triangle of $B_i$} if its vertices are all in $\bX_i$, or
\item an \emph{outer triangle of $B_i$} if none of its vertices is in $\bX_i$.
\end{itemize}
Notice that boundary triangles are the triangles that are cut by one of the grid lines. See Figure~\ref{fig:alpha_complex_intermediate} for an illustration of these types of triangles.

The union of all inner and boundary triangles of a zone $B_i$ and their faces is a simplicial subcomplex of $\DX$. 
We denote it by $K_i$.
A simplex from $\DX$ is either
\begin{itemize}
    \item a \emph{boundary simplex of $K_i$} if it is the face of a boundary triangle of $B_i$, 
    \item an \emph{inner simplex of $K_i$} if all of its vertices are lying in $\bX_i$ and it is not a boundary simplex of $K_i$, or 
    \item an \emph{outer simplex of $K_i$} if none of its vertices is in $\bX_i$ and it is not a boundary simplex of $K_i$.
\end{itemize}
In particular, notice that boundary, inner and outer triangles of $B_i$ are boundary, inner and outer simplices of $K_i$ respectively.

\begin{figure}
    \centering
    \begin{subfigure}[b]{0.45\textwidth}
        \centering
        \includegraphics[scale=0.5]{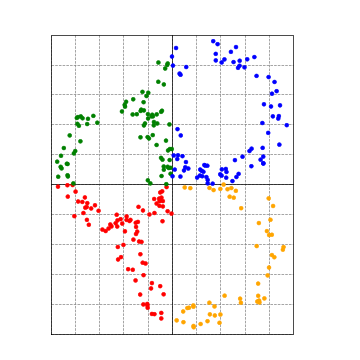}
        \caption{$\bX$ partitioned into four zones.}
        \label{fig:point_cloud_subdivided}
    \end{subfigure}  
    \hfill
    \begin{subfigure}[b]{0.45\textwidth}
        \centering
        \includegraphics[scale=0.3]{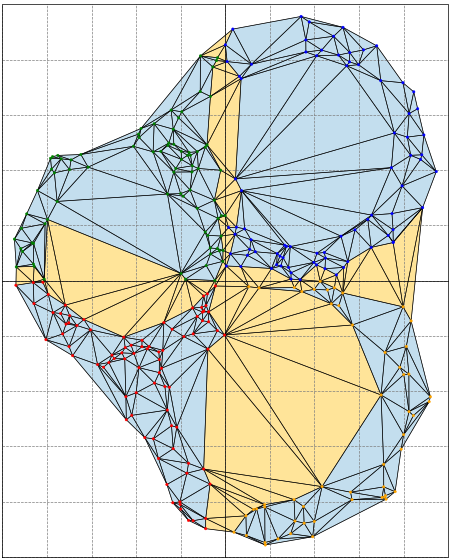}
        \caption{$\DX$ subdivided into four.}
        \label{fig:delaunay_triangulation_final}
    \end{subfigure}
    \caption{Division of $\bX$ into four zones and the induced cover for $\DX$. Intersections between subcomplexes are depicted in yellow in (b).}
    \label{fig:three graphs}
\end{figure}

The set of subcomplexes $\{K_i\}_{0\leq i < M}$ defines a cover of $\DX$, i.e. $\bigcup_{0\leq i < M} K_i = \DX.$
In general, fixed $i$, $K_i$ is not the Delaunay triangulation of the point cloud defined by the vertex set $V(K_i)$ as the union of simplices in $K_i$ might not be the convex hull of $V(K_i)$.
However, all triangles from $K_i$ satisfy the Delaunay condition. 
\begin{example} 
Recall the points $\bX$ and the Delaunay triangulation from Example~\ref{ex:alpha_filtered}. 
We take the rectangular subspace $B\subseteq \bR^2$ containing $\bX$ and divide it into four zones as shown in Subfigure~\ref{fig:point_cloud_subdivided}. 
This induces a cover of the Delaunay triangulation into four subcomplexes depicted in Figure~\ref{fig:delaunay_triangulation_final}.
\end{example}
Since the intersections of the covers play a crucial role in this article, we shortly explain how they can be characterised.
A simplex $\sigma$ lies in the intersection of subcomplexes $\bigcap_{j \in J} K_j$, where $J \subseteq \{0,\ldots, M-1\}$, if, and only if, it is a boundary simplex of all of them.
If there was one $u \in J$ such that $\sigma$ is not in the boundary of $K_{u}$, $\sigma$ would either be an outer simplex or an inner simplex of $K_{u}$. 
If it was an outer simplex it would not be a simplex of $K_{u}$. 
If it was an inner simplex it would not be a simplex of any $K_j$, $j\in J\backslash u$.
And thus, in both cases, $\sigma$ would not be in $\bigcap_{j \in J} K_j$.
Using this, we only have to look at the boundary simplices when computing the intersections.
\leavevmode

Next, we obtain further characterisations for higher intersections. 
First, notice that a  
triangle $\tau$
lies in an intersection $\bigcap_{j \in J} K_j$ if and only if $\tau$ has at least one vertex in $\bX_j$ for all $j \in J$. 
A direct consequence is that all intersections for which $\#J>3$ cannot contain triangles.  
Now, consider an edge $\gamma$ lying in an intersection $\bigcap_{j \in J} K_j$. 
It follows that $\gamma$ is a boundary simplex for all $K_i$ with $i\in J$. Equivalently, $\gamma$ is the face of a 
triangle with a vertex in $\bX_i$ for all $i \in J$. 
As $\gamma$ is an edge from $\DX$, there exist at most two triangles $\tau_1$ and $\tau_2$ sharing $\gamma$ as a common face.
Since both $\tau_1$ and $\tau_2$ involve $4$ different vertices, this implies that $\#J\leq 4$. 
Thus, intersections involving five or more cover elements do not contain edges.
\begin{lemma}\label{lem:dim-high-inters}
    If $\#J\geq 4$ then $\dim \big(\bigcap_{j \in J} K_j\big) \leq \max\{5-\#J, 0\}$. 
\end{lemma}

\subsection{Cover construction}\label{sub:alpha-covers}
In this subsection, we summarise the algorithm to obtain the covers of $\DX$ in parallel. Details can be found in~\cite{Lo2012}. At the end of the subsection we illustrate this part of the algorithm with the example point cloud first introduced in Section~\ref{sec:alpha_complex}.

As described in Subsection~\ref{sub:definition_cover}, the point cloud $\bX$ is subdivided into disjoint subsets $\bX_i$, $i=0,\ldots,M-1$, based on a rectilinear grid.
The number of subsets $M$ is 
such that each set $\bX_i$ is stored and handled by a different processor.
Details on possible choices for $M$ may be found in Section~\ref{sec:implementation}.

For each subset $\bX_i$ the corresponding Delaunay triangulation $\DXi{i}$ is computed.
Since the sets $\bX_i$ are disjoint and bounded by the bounding boxes $B_i$,
the complexes $D_i^0 \coloneqq \DXi{i}$ do not intersect with each other.
We therefore insert a layer of points from the neighbouring sets $\bX_j$ obtaining an updated Delaunay triangulation $D_i^1$.
The thickness of the layer depends on the underlying grid.
We call the resulting vertex set $\bX_i^1$ and the new, slightly bigger bounding box $B_i^1$.

Two complexes $D_i^1$ and $D_j^1$ intersect now if the original bounding boxes $B_i$ and $B_j$ have a common edge.
Furthermore, the complex $D_i^1$ consists of inner, boundary and possibly also outer triangles of $B_i$. 
We define the subcomplexes $K_i^1\subset D_i^1$ to be the subcomplexes consisting only of inner and boundary triangles of $B_i$ and their faces.
We have to check now whether these complexes are equal to the complexes $K_i$ from Subsection~\ref{sub:definition_cover}.
Assuming $\bX$ to be in general position, $K_i^1$ is equal to $K_i$, if and only if every triangle of $K_i^1$ is Delaunay with respect to $\bX$.
Using Delaunay's Lemma (Theorem \ref{thm:delaunay}), it suffices to check the Delaunay property for the boundary triangles.
These are Delaunay with respect to $\bX_i^1$ but their circumcircles might exceed $B_i^1$.
In that case, there might be other points in $\bX\backslash\bX_i^1$ lying in such a circumcircle making the triangle non-Delaunay.
Then another insertion of points into the triangulation is necessary.
We proceed in an iterative way. 
In the $k$-th round, we update the Delaunay triangulation $D_i^k$, the bounding box $B_i^k$, the vertex set $\bX_i^k$ and the resulting complex $K_i^k$.

Assume there is at least one boundary triangle in $K_i^1$ whose circumcircle exceeds $B_i^1$.
Starting with $k=1$, we do the following:
\begin{enumerate}
    \item Insert another layer of points surrounding $B_i^k$ into $D_i^k$,
    obtaining $B_i^{k+1}$, $\bX_i^{k+1}$ and $D_i^{k+1}$.
    \item Compute $K_i^{k+1}$ from $D_i^{k+1}$.
    \item Compute the circumcircles of all boundary triangles of $K_i^{k+1}$.
    \item Check whether there is at least one circumcircle exceeding $B_i^{k+1}$
    but intersecting $B \setminus B_i^{k+1}$.
    \item If yes, repeat. If not, $K_i^{k+1}=K_i$.
\end{enumerate}

For the insertion of new points communication between the processors is necessary.
The other steps may be handled individualy by the processors.

\begin{figure}[h]
    \begin{center}
    \includegraphics[width=0.7\textwidth]{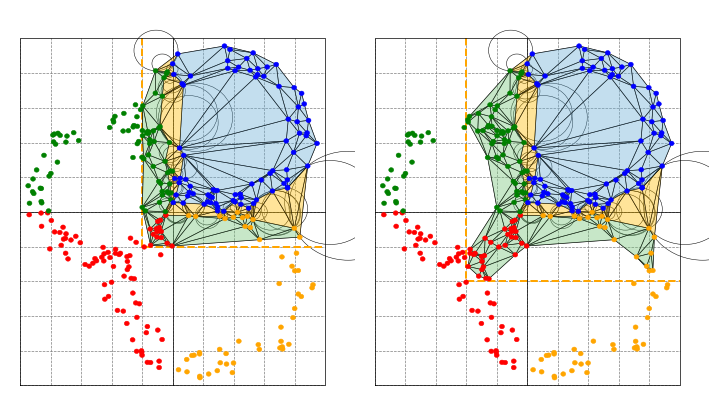}
    \caption{First and second step of iterative Delaunay Triangulation construction.
    Inner triangles are marked in blue, boundary triangles are marked in yellow and outer triangles are marked in green.}
    \label{fig:alpha_complex_intermediate}
    \end{center}
\end{figure}

\begin{example}\label{ex:cover_construction}
    Consider the point cloud introduced in Example~\ref{ex:alpha_filtered}.
    We want to subdivide the corresponding Delaunay triangulation into four subcomplexes.
    At first, we subdivide the point cloud into four disjoint sets based on a regular grid as depicted in Figure~\ref{fig:point_cloud_subdivided}.
    Then we compute the Delaunay triangulation of those.
    Next we iteratively add points to the triangulations and expand those until all boundary simplices satisfy the Delaunay condition.
    In Figure~\ref{fig:alpha_complex_intermediate} the first two steps of this iteration are depicted for the upper right subcomplex. 
    In this case, the iteration already stops after the second round.
    The circumcircles of the boundary simplices are depicted and after the second step they are all lying within the updated bounding box (dased yellow lines).
    The resulting four subcomplexes are depicted in Figure~\ref{fig:delaunay_triangulation_final}.
\end{example}
\leavevmode

We have now computed the subcomplexes $K_i$ for each zone $i$. 
For the Mayer-Vietoris spectral sequence we also need to compute the intersections of those, which will be explained in the next subsection.

\subsection{Intersection Computation}\label{sub:alpha-intersection}
In this section, we describe how each processor $i$ computes the intersections of the corresponding subcomplex $K_i$ with other subcomplexes. 
To do this, communication between the processors is necessary, but most of it may be computed in parallel.

Next, as mentioned in Subsection~\ref{sub:definition_cover}, given $j\neq i$, a simplex is in $K_i\cap K_j$ if and only if it is a boundary simplex for both $K_i$ and $K_j$.
Thus, for the computation of $K_i\cap K_j$ we only have to check the boundary simplices.
Furthermore, a triangle is in $K_i\cap K_j$ if and only if at least one of its vertices lies in $B_i$ and at least one of its vertices lies in $B_j$.
This makes assigning triangles and their faces to intersections very simple.

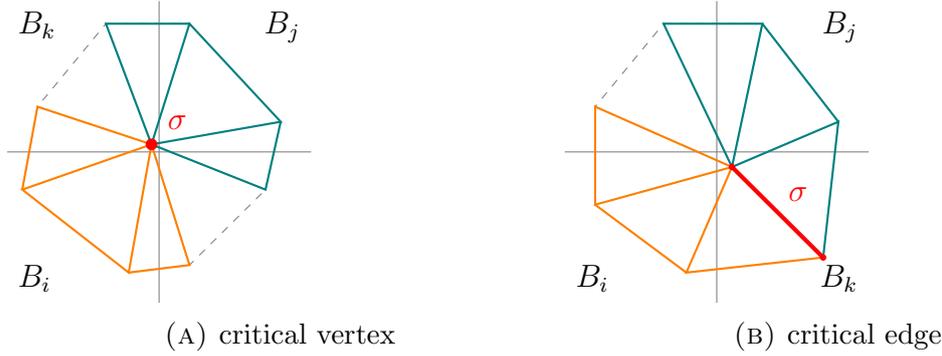
\begin{figure}
    \centering
    \begin{subfigure}{0.45\textwidth}
    \begin{tikzpicture}
        \coordinate (A) at (-0.1,0.1);
        \coordinate (B) at (0.4,1.7);
        \coordinate (C) at (1.6, 0.4);
        \coordinate (D) at (1.4, -0.5);
        \coordinate (E) at (0.4, -1.5);
        \coordinate (F) at (-0.4, -1.6);
        \coordinate (G) at (-1.8,-0.5);
        \coordinate (H) at (-1.6, 0.6);
        \coordinate (I) at (-0.7, 1.7);
        \draw[gray] (-2,0) -- (2,0);
        \draw[gray] (0,-2) -- (0,2);
        \draw (-2,-2) node [anchor=south west] {$B_i$};
        \draw (2,2) node [anchor=north east] {$B_j$};
        \draw (-2,2) node [anchor=north west] {$B_k$};
        \draw[gray, dashed] (A) -- (D) -- (E) -- (A);
        \draw[gray, dashed] (A) -- (H) -- (I) -- (A);
        \draw[gray, dashed] (D) -- (C);
        \draw[gray, dashed] (I) -- (B);
        \draw[gray, dashed] (G) -- (H);
        \draw[gray, dashed] (F) -- (E);
        \draw[color=teal, thick] (A) -- (B) -- (C) -- (A) -- (D) -- (C);
        \draw[color=teal, thick] (A) -- (I) -- (B);
        \draw[color=orange, thick] (A) -- (F) -- (G) -- (A) -- (H) -- (G);
        \draw[color=orange, thick] (A) -- (E) -- (F);
        \filldraw[red] (A) circle (2pt);
        \draw ($(A)!0.2!(C)$) node [text=red, anchor=south] {$\sigma$};
        \end{tikzpicture}
        \caption{critical vertex}
        \end{subfigure}
        \begin{subfigure}{0.45\textwidth}
        \begin{tikzpicture}
        \coordinate (A) at (0.2,-0.2);
        \coordinate (B) at (0.6,1.7);
        \coordinate (C) at (1.6, 0.4);
        \coordinate (D) at (1.4, -1.4);
        \coordinate (F) at (-0.4, -1.6);
        \coordinate (G) at (-1.6,-0.7);
        \coordinate (H) at (-1.6, 0.6);
        \coordinate (I) at (-0.7, 1.7);
        \draw[gray] (-2,0) -- (2,0);
        \draw[gray] (0,-2) -- (0,2);
        \draw (-2,-2) node [anchor=south west] {$B_i$};
        \draw (2,2) node [anchor=north east] {$B_j$};
        \draw (2,-2) node [anchor=south east] {$B_k$};
        \draw[gray, dashed] (A) -- (H) -- (I) -- (A);
        \draw[gray, dashed] (I) -- (B);
        \draw[gray, dashed] (G) -- (H);
        \draw[gray, dashed] (F) -- (G) -- (A);
        \draw[gray, dashed] (C) -- (B) -- (A);
        \draw[color=orange, thick] (D) -- (F) -- (A) -- (G) -- (H) -- (A);
        \draw[color=orange, thick] (F) -- (G);
        \draw[color=teal, thick] (D) -- (C) -- (A) -- (I) -- (B) -- (A);
        \draw[color=teal, thick] (B) -- (C);
        \filldraw[red] (A)  circle (1pt);
        \filldraw[red] (D)  circle (1pt);
        \draw[red, ultra thick] (A) -- (D);
        \draw ($(A)!0.5!(D)$) node [text=red, anchor=south west] {$\sigma$};
        \end{tikzpicture}
        \caption{critical edge}
        \end{subfigure}
    \caption{Two examples of shared simplices between $K_i$ (orange) and $K_j$ (dark green) which aren't faces of shared 2-simplices}
    \label{fig:critical_edge_vertex}
\end{figure}

Unfortunately, there may be vertices and edges in intersections that are not faces of a shared triangle. 
Two such cases are illustrated in Figure~\ref{fig:critical_edge_vertex}.
In the first example, the vertex $\sigma$ lies in both $K_i$ and $K_j$ even though the two complexes do not share a triangle. 
In the second example, different boundary triangles of $K_i$ and $K_j$ touch in the edge $\sigma$. Again, $K_i$ and $K_j$ do not share any boundary triangles.

We want to characterise these vertices and edges in more detail.
Let $\sigma$ be such a simplex in the intersection of subcomplexes $K_i$ and $K_j$, $i\neq j$.
For $K_i$, we know that $\sigma$ is a boundary simplex, and thus, there is a boundary triangle $\tau$ of which $\sigma$ is a face.
However, since $\tau$ is not shared with $K_j$, there must be at least another subcomplex $K_k$ that shares $\tau$ with $K_i$.
From the above, we can also deduce that $\sigma$ does not have vertices in either $K_i$ or $K_j$. 
This observation leads to the following definition.

\begin{definition}[Critical Vertex and Critical Edge]
    Fix $i,j\in\{0,\ldots,M-1\}$, $i\neq j$.
    A \textit{critical vertex} of $K_i\cap K_j$ is defined to be a vertex of $K_i\cap K_j$ not lying in either $\bX_i$ or $\bX_j$.
    Similarly, a \textit{critical edge} of $K_i\cap K_j$ is defined to be an edge of $K_i\cap K_j$ whose vertices are neither in $\bX_i$ nor in $\bX_j$. 
\end{definition}

Now, let $\sigma$ be a critical vertex of $K_i\cap K_j$ which is a vertex from $\bX_k$.
From this, we can deduce that $\sigma$ is the face of a triangle in both $K_k\cap K_i$ and the face of another triangle in $K_k \cap K_j$, and thus $\sigma \in K_k\cap K_i \cap K_j$.

\leavevmode

With these observations in mind, processor $i$ computes  intersections $K_i \cap K_j$ and  $K_i\cap K_j \cap K_k$ in five steps:
\begin{enumerate}

\item Compute all shared triangles and their faces by going through the list of boundary simplices and adding them to a two-fold intersection $K_i \cap K_j$ whenever at least one vertex is in the other point set $\bX_j$. 
In the same step, add a triangle and its faces to $K_i\cap K_j \cap K_k$ whenever one vertex is in $\bX_j$ and another one in $\bX_k$.
Notice that there is always at least one vertex in $\bX_i$.

\item Compare two-fold intersections to detect shared edges and vertices in three-fold intersections that are not faces of triangles in these three-fold intersections.
Add the detected simplices to the corresponding three-fold intersections.
In Figure~\ref{fig:critical_edge_vertex}, in both cases, the processor $k$ detects the respective $\sigma$ to be in the three-fold intersection $K_i\cap K_j \cap K_k$.

\item Share the preliminary three-fold intersections between all processors.
In Figure~\ref{fig:critical_edge_vertex}, processor $k$ sends the information of $\sigma$ being a shared simplex of $K_i$, $K_j$ and $K_k$ to processors $i$ and $j$.

\item Update the three-fold intersections and subsequently the two-fold intersections by adding the missing simplices.
In Figure~\ref{fig:critical_edge_vertex}, processor $i$ adds the respective $\sigma$ to the intersection $K_i\cap K_j \cap K_k$ and subsequently to the intersection $K_i \cap K_j$.
\end{enumerate}

For this work we do not need to compute higher intersections.

\subsection{Computation of Filtration Values}\label{sec:alpha-fv}
So far, we have constructed a cover for $\DX$. 
We obtain a filtered cover for $\DX$ by intersecting each subcomplex $K_i$ with the filtered complex $\AXn$,
obtaining subcomplexes $\rA_i^r \coloneqq K_i \cap \AX{r}$ for $r=0,\ldots, N$, and we denote as $\rA_i$ the resulting filtration on $K_i$.

Recall that the filtration value of any triangle and any Gabriel edge is given by the squared radius of its circumcircle.
The filtration value of a non-Gabriel edge is given by the minimum filtration values of its cofaces.
And at last, the filtration value of a vertex is zero.
Thus, computing the filtration values of all triangles can be done in parallel as it is a local property.
The same holds true for all edges whose cofaces are contained in $K_i$.

However, we need to be a bit more careful with the non-Gabriel edges that are on the boundary of a subcomplex $K_i$ and with those in the intersections of two or more subcomplexes.
In both cases, the processor computing the filtration value for the particular edge might not be able to detect that it is a non-Gabriel edge. 
That is, for an edge to be Gabriel, it suffices to check whether the vertices of the adjacent triangles (the cofaces of the edge) that are opposite of the edge are lying outside of its circumcircle.
But, if an edge is on the (topological) boundary of a subcomplex $K_i$, the particular vertex making it non-Gabriel might be the vertex of the triangle not lying in $K_i$. If this is the case, such non-Gabriel edge must be in a subcomplex $K_j$ such that it also contains the vertex that lies in its circumcircle; we illustrate this with an example.

\begin{example}
Consider an intersection $K_i\cap K_j$ containing a non-Gabriel edge, depicted in red in Figure~\ref{fig:non_gabriel_shared_edge}.
In such a case, processor $i$, which is computing the filtration values for $\rA_i$, cannot detect that the edge is actually non-Gabriel.
\begin{figure}[ht]
    \centering
    \begin{tikzpicture}
        \coordinate (A) at (0,3);
        \coordinate (B) at (0,-1);
        \coordinate (C) at (1.4,1);
        \coordinate (D) at (-3.4,1.6);
        \coordinate (E) at ($(A)!0.5!(B)$);
        \draw[gray, dashed] (-1,-1.2) -- (-1,3.2);
        \draw[teal,dashed] ($(A)!0.5!(B)$) circle (2.0);
        \path let \p1=($(A)-(B)$),\p2=($(A)!0.5!(B)$),
            \p3=($(A)-(C)$),\p4=($(A)!0.5!(C)$), 
            \n1={(-(\x2*\x3)+\x3*\x4+\y3*(-\y2 +\y4))/(\x3*\y1-\x1*\y3)} in coordinate (F) at ({\x2+\n1*\y1},{\y2-\n1*\x1});
        \path let \p1=($(A)-(B)$),\p2=($(A)!0.5!(B)$),
            \p3=($(A)-(D)$),\p4=($(A)!0.5!(D)$), 
            \n1={(-(\x2*\x3)+\x3*\x4+\y3*(-\y2 +\y4))/(\x3*\y1-\x1*\y3)} in coordinate (G) at ({\x2+\n1*\y1},{\y2-\n1*\x1});
        \filldraw[orange, opacity=.3] (A) -- (B) -- (D) -- (A);
        \filldraw[orange, opacity=.3] (A) -- (-0.3,3.2) -- (-3.6,3.2) -- (D) -- (A);
        \filldraw[orange, opacity=.3] (B) -- (-0.5,-1.2) -- (-3.2,-1.2) -- (D) -- (B);
        \draw (-2.5,1.3) node [text=orange, anchor=west] {$K_i\cap K_j$};
        \draw (-4.5,0) node [anchor=south west] {$K_i$};
        \draw (1.5,0) node [anchor=south west] {$K_j$};
        \coordinate (auxAC) at ($(A)!0.5!(C)$);
        \coordinate (auxBC) at ($(B)!0.5!(C)$);
        \coordinate (auxAD) at ($(A)!0.5!(D)$);
        \coordinate (auxBD) at ($(B)!0.5!(D)$);
        \draw[thick] (A) -- (C) -- (B) -- (D) -- (A);
        \draw[thick,red] (A) -- (B);
        \foreach \i in {A,B,C,D} {
            \filldraw (\i) circle (2pt);
        }
        \draw[red] (0.1,2) -- (1,2.5);
        \draw (1,2.5) node [text=red,anchor=west] {critical non-Gabriel edge};
        \draw[red] (1.45,1.2) -- (1.55,1.45);
        \draw (1.4,1.6) node [text=red,anchor=west] {vertex in its circumcircle};
        \draw[red,thick] (C) circle (4pt);
        \draw[thick, dotted] (D) -- (-3.6,3.2);
        \draw[thick, dotted] (D) -- (-3.2,-1.2);
        \draw[thick, dotted] (D) -- (-4,2);
        \draw[thick, dotted] (D) -- (-4,0.9);
        \draw[thick, dotted] (A) -- (-0.3,3.2);
        \draw[thick, dotted] (A) -- (0.2,3.2);
        \draw[thick, dotted] (B) -- (-0.5,-1.2);
        \draw[thick, dotted] (B) -- (0.3,-1.2);
        
    \end{tikzpicture}
    \caption{Critical non-Gabriel edge (red). The edge cannot be identified as non-Gabriel by $K_i$ or $K_i\cap K_j$ since the vertex causing it to be non-Gabriel does not lie in $K_i$ (and therefore also not in $K_i\cap K_j$).}
    \label{fig:non_gabriel_shared_edge}
\end{figure}
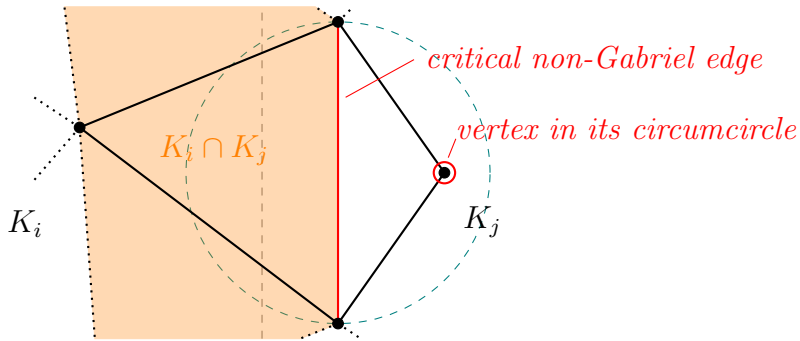
\end{example}

In order to solve this issue, while still computing the filtration values in parallel as much as possible, process $i$ proceeds as follows:
\begin{enumerate}
\item Compute the filtration values just for the complex $\rA_i$.
In this first step, the process stores a list identifying non-Gabriel edges from $\rA_i$. 
Notice that, at this point, $\rA_i$ might be inconsistent with the filtration computed by other processes.

\item Compute the filtration values for pairwise intersections and, whenever dealing with an edge, check whether it is in the list of non-Gabriel edges.
If an edge is found in the list, check whether it is also non-Gabriel in the intersection. If this check fails, assign the correct filtration value and identify it as a \emph{critical non-Gabriel edge}.

\item Exchange these lists of critical non-Gabriel edges between the processors and update the filtration values, if necessary.
\end{enumerate}
In that way, we assign the same filtration value to every simplex as in $\AXn$.

\begin{example}
In Figure~\ref{fig:non_gabriel_shared_edge}, when only looking at $K_i\cap K_j$, processor $j$ would not be able to identify the non-Gabriel edge as such.
In Figure~\ref{fig:non_gabriel_shared_edge}, processor $j$ shares its list with processor $i$ which subsequently updates the filtration value assigned to the red edge in $K_i$ and in all intersections that contain it.
\end{example}

\section{PerMaViss on Alpha-Grid covers}\label{sec:permaviss_alpha}
In this Section, we consider the persistence Mayer-Vietoris spectral sequence from~\cite{Torras2023} in the particular case of alpha complex covers introduced in Section~\ref{sec:alpha-construction}.

\subsection{Double complex and total complex}~\label{sub:double-total}
Let $\AXn$ be an alpha complex covered by $\cU=\{\rA_i\}_{i=0}^{M-1}$, as done in Section~\ref{sec:alpha-construction}. Before delving into details, we denote by $N(p)$ the cardinality of $\cN(\cU)^p$ for all $p\geq 0$. Also, we organise the simplices from $\cN(\cU)^p$ into sequences $(\sigma_1, \sigma_2, \ldots, \sigma_{N(p)})$ for ease. We will work with the tensor product $\otimes_{\bZ_2}$ on $\bZ_2$-vector spaces and, for ease, we write $\otimes$ instead of $\otimes_{\bZ_2}$.
We consider the \emph{double complex} whose terms are given by the following subspaces of $C_p(\cN(\cU))\otimes C_q(\DX)$, for $p,q \geq 0$,
\[
C_{p,q}(\AXn, \cU) = \bigoplus_{i=1}^{N(p)} \sigma_i \otimes C_q(\rA(\sigma_i))\, \mbox{, where }
\rA(\sigma_i) = \bigcap_{j \in \sigma_i}\rA_j\,.
\]
Now, since $\AXn$ is of dimension $2$, the same follows for all subcomplexes $\rA(\sigma)$ for all $\sigma \in \cN(\cU)$. In particular, $C_{p,q}(\AXn, \cU)=0$ for all $q\geq 3$. In addition, by Lemma~\ref{lem:dim-high-inters}, we also have that $C_{p,q}(\AXn, \cU)=0$ for all pairs $(p,q)$ such that both $q\geq 1$ and $p+q\geq 5$ hold. Denoting the dimension of $\cN(\cU)$ by $D$ (normally $D\leq 3$), we also have $C_{p,q}(\AXn, \cU)=0$ for all $p\geq D+1$. 
For ease, we write $C_{p,q}$ instead of $C_{p,q}(\AXn, \cU)$ from now on.
So we might depict the nontrivial part of this double complex as in the following commutative diagram

\begin{equation}\label{diag:double-cpx}
\begin{tikzcd}[column sep=0.7cm, /tikz/column 7/.style={column sep=-0.5em}]
    C_{0,2} \ar[d, "d^V"] &
    C_{1,2} \ar[l, "d^H", swap] \ar[d, "d^V"] & 
    C_{2,2} \ar[l, "d^H", swap] \ar[d, "d^V"] \\
    C_{0,1} \ar[d, "d^V"] &
    C_{1,1} \ar[l, "d^H", swap] \ar[d, "d^V"] & 
    C_{2,1} \ar[l, "d^H", swap] \ar[d, "d^V"] & 
    C_{3,1} \ar[l, "d^H", swap] \ar[d, "d^V"] \\
    C_{0,0} &
    C_{1,0} \ar[l, "d^H", swap] & 
    C_{2,0} \ar[l, "d^H", swap] & 
    C_{3,0} \ar[l, "d^H", swap] &
    C_{4,0} \ar[l, "d^H", swap] & 
    \cdots \ar[l, "d^H", swap] & 
    C_{D,0} \ar[l, "d^H", swap] & ,
\end{tikzcd}
\end{equation}
where the arrows are given by the \emph{vertical differentials} $d^V_{p,q}:C_{p,q}\rightarrow C_{p,q-1}$ and the \emph{horizontal differentials} $d^V_{p,q}:C_{p,q}\rightarrow C_{p-1,q}$ respectively; we will often omit the subindices as done in Diagram~(\ref{diag:double-cpx}). Next, we will describe these differentials with detail. 

\begin{figure}
    \centering
    \begin{minipage}{.3\textwidth}
        \centering 
        \begin{subfigure}[b]{\textwidth}
            \centering
            \includegraphics[width=0.6\textwidth]{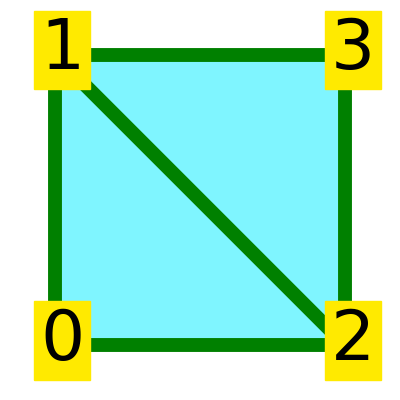}
            \caption{Nerve $\cN(\cU)$}
            \label{subfig:nerve-grid}
        \end{subfigure}
        \vfill
        \begin{subfigure}[b]{\textwidth}
            \centering
            \includegraphics[width=\textwidth]{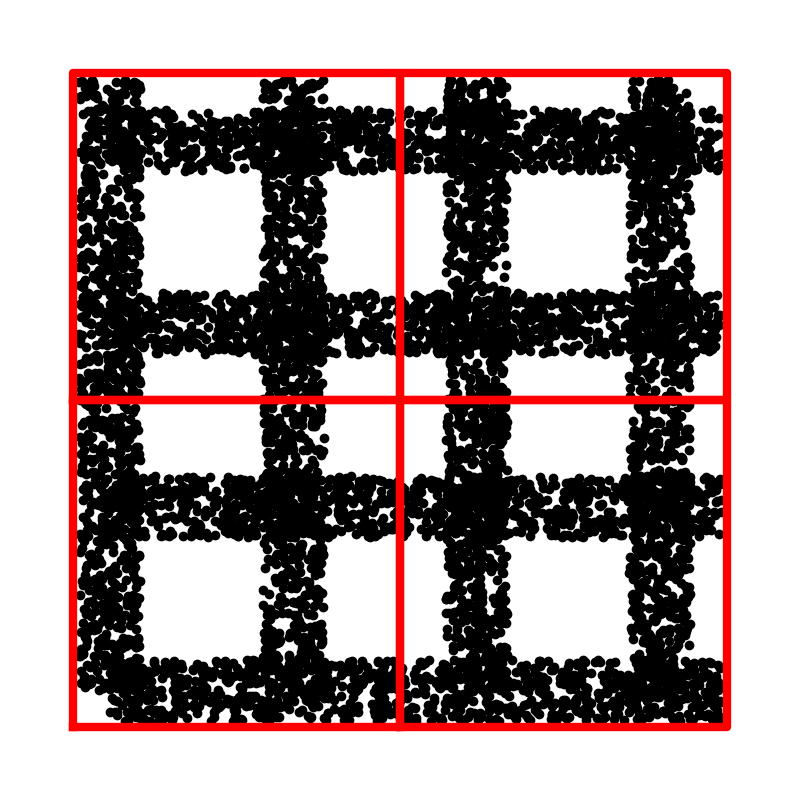}
            \caption{Grid points and zones}
            \label{subfig:zones-grid}
        \end{subfigure}
    \end{minipage}
    \hfill
    \begin{minipage}{.68\textwidth}
    \begin{subfigure}[b]{\textwidth}
        \centering
        \includegraphics[width=\textwidth]{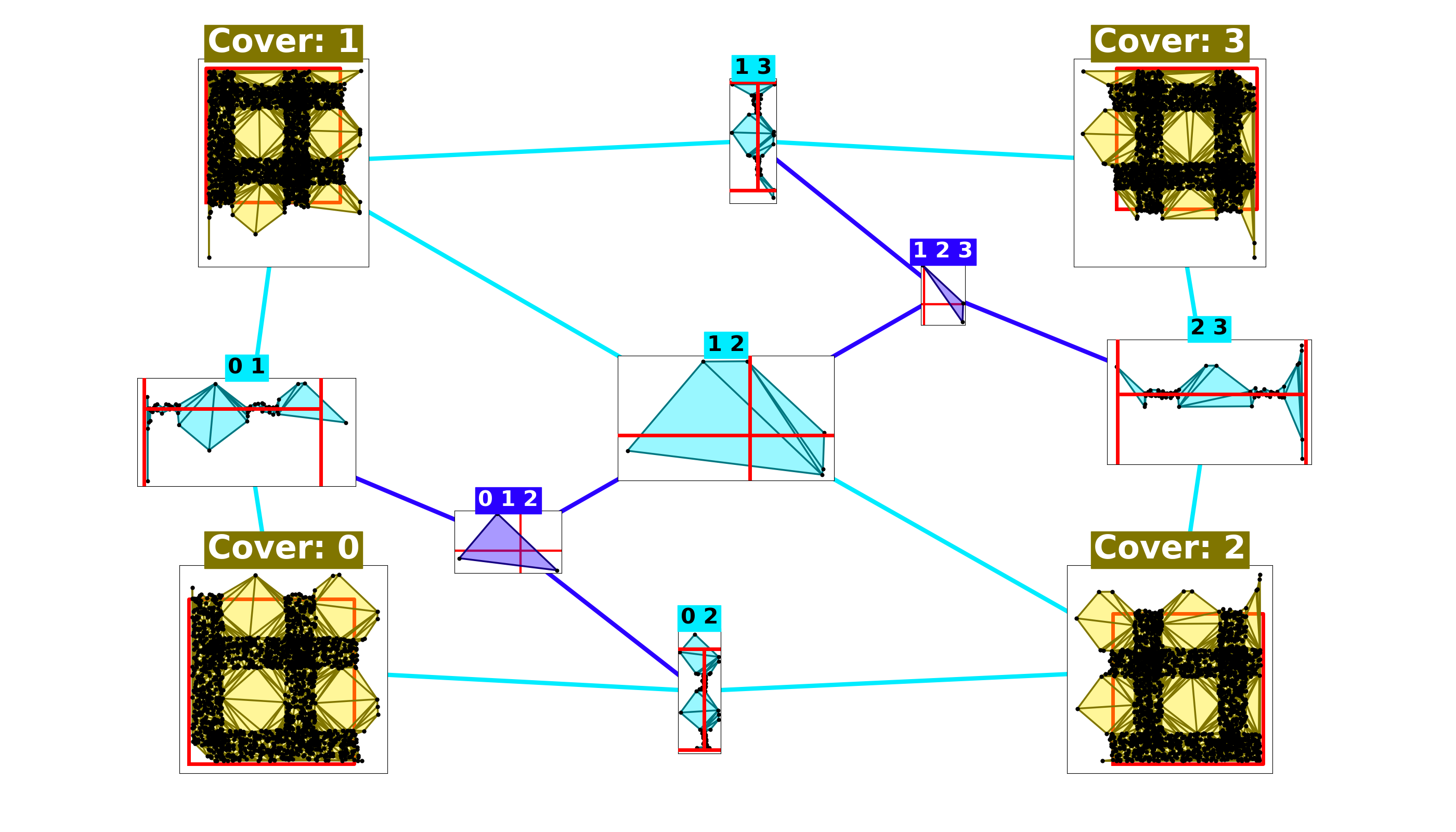}
        \caption{Covering regions with inclusion indicated by lines.}
        \label{sub:grid-cover-sets}
    \end{subfigure}
    \end{minipage}
    \caption{Point cloud divided into four zones and its nerve.}\label{fig:grid-cover-example}
\end{figure}

The vertical differentials are linear morphisms such that, given $\sigma \otimes a \in C_{p,q
}$, satisfy $d^V_{p,q}(\sigma \otimes a) = \sigma \otimes d_{\sigma, q}(a)$, where $d_{\sigma, q}:C_q(K(\sigma))\rightarrow C_{q-1}(K(\sigma))$ are the $q$-differentials on the complex $K(\sigma)$; we denote their associated matrices as $D_{\sigma, q}$.
Thus, the matrix associated to $d^V_{p,q}$ is given by the diagonal block matrix
\[
D^V_{p,q} = \left[
\begin{array}{*{5}{w{c}{0.9cm}}}
    \blCell D_{\sigma_1, q} & 0 & 0 & \cdots & 0 \\
    0 & \blCell D_{\sigma_2, q} & 0 & \cdots & 0 \\
    0 & 0 & \blCell \ddots &  \cdots & 0 \\
    \vdots & \vdots & \vdots & \blCell \ddots & \vdots \\
    0 & 0 & 0 & \hdots & \blCell D_{\sigma_N, q}
\end{array}
\right]\,.
\]
The horizontal differentials are linear morphisms such that, given $\sigma \otimes a \in C_{p,q
}$, satisfy $d^H_{p,q
}(\sigma \otimes a) = d^{\cN(\cU)}_p(\sigma)\otimes a$.
The horizontal differentials are handled by a block matrix $D^H_{p,q}$ which depends on the nerve $\cN(\cU)$.
Essentially, given any pair of simplices $\tau$ and $\sigma$ from $\cN(\cU)$ and such that $\vert\sigma\vert = \vert\tau\vert+1$, we denote by $D^\tau_{\sigma, q}$ the block from $D^H_{p,q}$ associated to the columns from the summand $\sigma$ and the rows from the summand $\tau$.
Now, if $\tau \preceq \sigma$, then $D^\tau_{\sigma, q}$ is the block associated to  $S_q(K(\sigma))\rightarrow S_q(K(\tau))$ which is induced by the inclusion $K(\sigma) \hookrightarrow K(\tau)$. 
On the other hand, if $\tau \not\preceq \sigma$, then $D^\tau_{\sigma, q}=0$.
\begin{example}
    Consider the point cloud $\bX$ as well as the zone division depicted in Subfigure~(\ref{subfig:zones-grid}). In this case, the cover $\cU$ only has two nontrivial triple intersections, which leads to the nerve $\cN(\cU)$ depicted in Subfigure~(\ref{subfig:nerve-grid}). The four covering regions and their intersections are depicted in Subfigure~(\ref{sub:grid-cover-sets}). We organise the nerve simplices into sequences following the lexicographic order; e.g. nerve edges are in a sequence $([0,1], [0,2], [1,2], [1,3], [2,3])$. Then, the matrices associated to the first and second horizontal differentials have the following block forms
    \[
  D^H_{1,q} = \left[
\begin{array}{ccccc}
    \blCell D^{[0]}_{[0,1], q} & \blCell D^{[0]}_{[0,2], q} & 0 & 0 & 0 
    \\
    \blCell D^{[1]}_{[0,1], q} & 0 & \blCell D^{[1]}_{[1,2], q} & \blCell D^{[1]}_{[1,3], q} & 0 
    \\
    0 & \blCell D^{[2]}_{[0,2], q} & \blCell D^{[2]}_{[1,2], q} &  0 & \blCell D^{[2]}_{[2,3], q} 
    \\
    0 & 0 & 0 & \blCell D^{[3]}_{[1,3], q} & \blCell D^{[3]}_{[2,3], q}
\end{array}
\right]\,,  
D^H_{2,q} = \left[
\begin{array}{ccc}
    \blCell D^{[0,1]}_{[0,1,2], q} & 0
    \\
    \blCell D^{[0,2]}_{[0,1,2], q} & 0
    \\
    \blCell D^{[1,2]}_{[0,1,2], q} & \blCell D^{[1,2]}_{[1,2,3], q}
    \\
    0 & \blCell D^{[1,3]}_{[1,2,3], q}
    \\
    0 & \blCell D^{[2,3]}_{[1,2,3], q}
\end{array}
\right]\,.
    \]
\end{example}
\leavevmode
One can check that $d^V \circ d^V=0$ and $d^H \circ d^H=0$ and $d^H \circ d^V+d^V \circ d^H=0$, which is equivalent to saying that $C_{p,q}(\AXn, \cU)$ is a double complex.

Since $C_{p,q}(\AXn, \cU)$ is a double complex, we can define the \emph{total complex} which is given by the terms, for all $m\geq 0$,
\[
C^\Tot_m(\AXn, \cU) = 
\bigoplus_{p+q=m} C_{p,q}(\AXn, \cU)
\]
and differentials 
$
d^\Tot_m=\bigoplus_{p+q=m} (d^V_{p,q} + d^H_{p,q})
$.
From the properties of $d^V$ and $d^H$, one can check that $d^\Tot$ is indeed a differential; i.e. $d^\Tot \circ d^\Tot=0$. Thus, one might compute the homology of the total complex and, as explained in~\cite{Torras2023}, we have that $\Ho_m(C^\Tot(\AXn, \cU))\simeq \PH_m(\AXn)$ for all $m\geq 0$. Thus, instead of directly computing the barcode of $\PH_m(\AXn)$, we can obtain the same result by computing the homology of the total complex. Fortunately, the Mayer-Vietoris spectral sequence gives a procedure that executes this computation by breaking it down into simple and computable steps.

\subsection{The persistence Mayer-Vietoris spectral sequence}~\label{sub:particular_permaviss}
To efficiently compute $\Ho_m(C^\Tot(\AXn, \cU))$, we consider the Mayer--Vietoris spectral sequence $E^n_{p,q}(\AXn, \cU)$, with $n\geq 0$, which is an algebraic device that breaks down the computation of $\Ho_m(C^\Tot(\AXn, \cU))$ into blocks that are deduced from the cover $\cU$. More precisely, 
 by computing the Mayer-Vietoris spectral sequence $E^r_{p,q}(\AXn, \cU)$ we obtain the persistent homology of the total complex $C^\Tot(\AXn, \cU)$ which is isomorphic to the persistent homology of $\AXn$. 
The zero page is given by the terms $E^0_{p,q}(\AXn, \cU) = C_{p,q}$. For ease, we denote terms as $E^r_{p,q}$ for all $r\geq 0$. 
The zero page differentials $d^0_{p,q}:E^0_{p,q}\rightarrow E^0_{p,q-1}$ are given by the aforementioned vertical differentials $d^V_{p,q}$ or zero if they involve terms outside Diagram~(\ref{diag:double-cpx}).
Computing homology with respect to these zero page differentials amounts to computing persistent homology independently over all complexes $\rA(\sigma)$ with $\sigma \in \cN(\cU)$. Thus, the first page terms are
\[
E^1_{p,q}(\AXn, \cU) = \bigoplus_{i=1}^{N(p)} \sigma_i \otimes \PH_q\big( \rA(\sigma_i)\big)\,.
\]
Now, recall that $\DX$ is acyclic and it is embedded into $\bR^2$. Hence, an elementary observation is that, for any subcomplex $L \subseteq \DX$, we have that $\Ho_2(L)=0$. 
Therefore, for any filtered subcomplex $L_*\subseteq \AXn$,  $\PH_k(L_*)=0$ for all $k\geq 2$.
In particular, $E^1_{p,q}=0$ for all $q\geq 2$.  This implies that, recalling the Double complex~(\ref{diag:double-cpx}), the nontrivial first page terms are as follows
\begin{equation}\label{diag:first_page}
\begin{tikzcd}[column sep=0.7cm, row sep=0.3cm, /tikz/column 7/.style={column sep=-0.5em}]
    E^1_{0,1} &
    E^1_{1,1} \ar[l, "d^1_{1,1}", swap] & 
    E^1_{2,1} \ar[l, "d^1_{2,1}", swap] & 
    E^1_{3,1} \ar[l, "d^1_{3,1}", swap] 
    \\
    E^1_{0,0} &
    E^1_{1,0} \ar[l, "d^1_{1,0}", swap] & 
    E^1_{2,0} \ar[l, "d^1_{2,0}", swap] & 
    E^1_{3,0} \ar[l, "d^1_{3,0}", swap] &
    E^1_{4,0} \ar[l, "d^1_{4,0}", swap] & 
    \cdots \ar[l, "d^1_{5,0}", swap] & 
    E^1_{D,0} \ar[l, "d^1_{D,0}", swap] & ,
\end{tikzcd}
\end{equation}
where $d^1_{p,q}:E^1_{p,q}\rightarrow E^1_{p-1,q}$ denote the first page differentials. These differentials are induced by the horizontal differentials $d^H_{p,q}$ and,
since the direct sum structure from the zero page is carried on to the first page, their matrices follow the same block matrix structure as $D^H_{p,q}$. To state this more precisely, calling $D^1_{p,q}$ the matrix associated to $d^1_{p,q}$, we also denote by $D^1_{q,\tau, \sigma}$ the block from $D^1_{p,q}$ associated to the columns from the summand $\sigma$ and the rows from the summand $\tau$. In particular, if $\tau \preceq \sigma$ then $D^1_{q,\tau, \sigma}$ is the matrix associated to $\PH_q(K(\sigma))\rightarrow \PH_q(K(\tau))$ and if $\tau \not\preceq \sigma$ then $D^1_{q,\tau,\sigma}$ is zero. 
In Appendix~\ref{app:associated-matrix} we outline a procedure for obtaining such matrices.

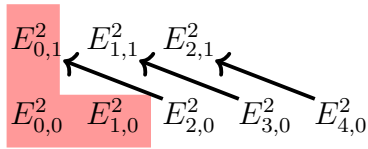
\begin{figure}[h]
\centering
\begin{tikzpicture}
  \foreach \i in {0,1,...,4}
  {
    \foreach \j in {0,1}
    {
      \node (v\i_\j) at (\i,\j) {};
    }
  }
  \fill[red!40] (-0.4,1.5)--(0.3,1.5)--(0.3,0.3)--(1.5,0.3)--(1.5,-0.4)--(-0.4,-0.4)--cycle;
  \foreach \i in {0,1,...,4}
  {
    \foreach \j in {0,1}
    {   \ifnum \i>4
            \ifnum \i=5
                \node at (v\i_\j) {0};
            \else 
                \node at (v\i_\j) {$\cdots$};
            \fi
        \else
            \pgfmathparse{(\i>2)&&(\j>0)}
            \ifnum \pgfmathresult=1
                \node at (v\i_\j) {};
            \else 
                \node at (v\i_\j) {$E^2_{\i,\j}$};
            \fi
        \fi
    }
  }
  \draw[line width=0.5mm, ->] ($(v2_0.north west)+(-0.2,0.1)$) -- ($(v0_1.south east)+(0.2,-0.1)$);
  \draw[line width=0.5mm, ->] ($(v3_0.north west)+(-0.2,0.1)$) -- ($(v1_1.south east)+(0.2,-0.1)$);
  \draw[line width=0.5mm, ->] ($(v4_0.north west)+(-0.2,0.1)$) -- ($(v2_1.south east)+(0.2,-0.1)$);
  
\end{tikzpicture}
\caption{Second page of the persistence Mayer-Vietoris spectral sequence. }
\label{fig:third_page_collapse}
\end{figure}

Computing homology with respect to the first page differentials, we obtain the second page terms $E^2_{p,q}=\Ker(d^1_{p,q})\big/ \Ima(d^1_{p+1,q})$. Now, recalling the first page structure from Diagram~(\ref{diag:first_page}), the nontrivial second page terms are distributed as shown in Figure~\ref{fig:third_page_collapse}. There are a few things to explain. First, the arrows denote the second page differentials $d^2_{p,q}:E^2_{p,q}\rightarrow E^2_{p-2,q+1}$ which might be nontrivial. Such differentials are trivial if $C_{p-1,q+1}=0$ and so, recalling Diagram~(\ref{diag:double-cpx}), $d^2_{p,q}=0$ for $p\geq 5$; i.e. there are only three possible nontrivial second page differentials. Further, computing homology with respect to these, we reach the third page, where all arrows must vanish and, since the same happens with higher pages, the spectral sequence \emph{collapses} on the third page. Next, since $\PH_q(\AXn)=0$ for all $q\geq 0$, the only possible nontrivial third page terms are in positions $(0,0)$, $(1,0)$ and $(0,1)$, which are within the red ``L" shaped region depicted in Figure~\ref{fig:third_page_collapse}.

By our deductions, it also follows that $E^2_{3,0} \simeq E^2_{1,1}$ while $E^2_{2,1} \simeq E^2_{4,0}$ and $\Ker(d_{2,0}^2)=0$. In the next subsection we will show conditions which guarantee that all second page differentials vanish, and so $E^2_{p,q}$ collapses on the second page, leaving only the terms $(0,0)$, $(1,0)$ and $(0,1)$. Also, in Example~\ref{ex:second-page-differential} we will see that this is not always the case. 

\subsection{Second Page Collapse}\label{sub:second_collapse}
\begin{figure}
    \centering
    \centering
    \begin{subfigure}[b]{0.25\textwidth}
    \centering 
    \includegraphics[width=\textwidth]{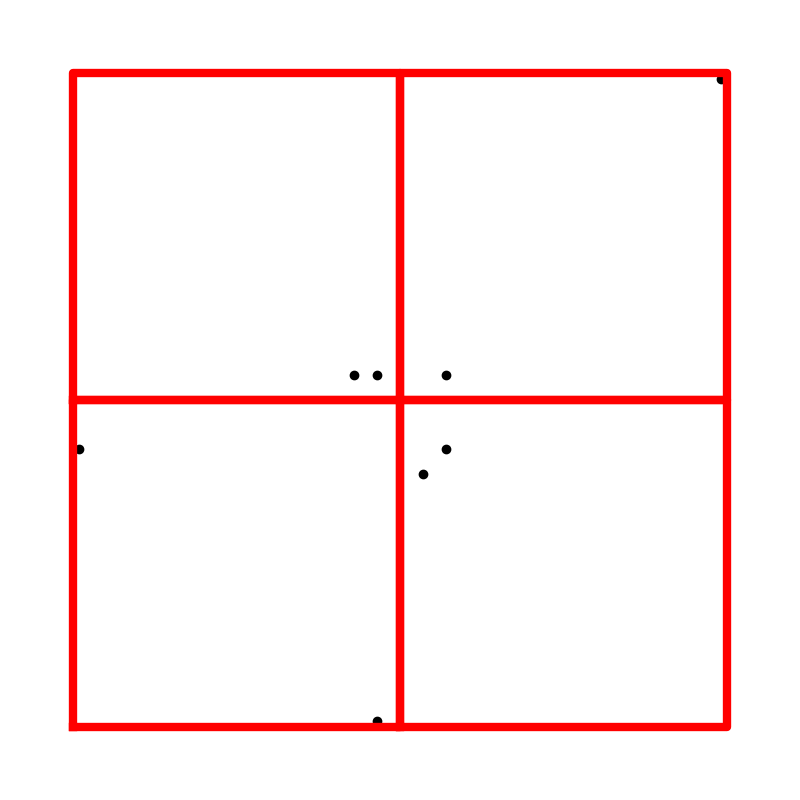}
    \caption{For zones covering $8$ points.}
    \label{subfig:zones-second-page}
    \end{subfigure}
    \hfill
    \begin{subfigure}[b]{0.7\textwidth}
    \centering \includegraphics[width=\textwidth]{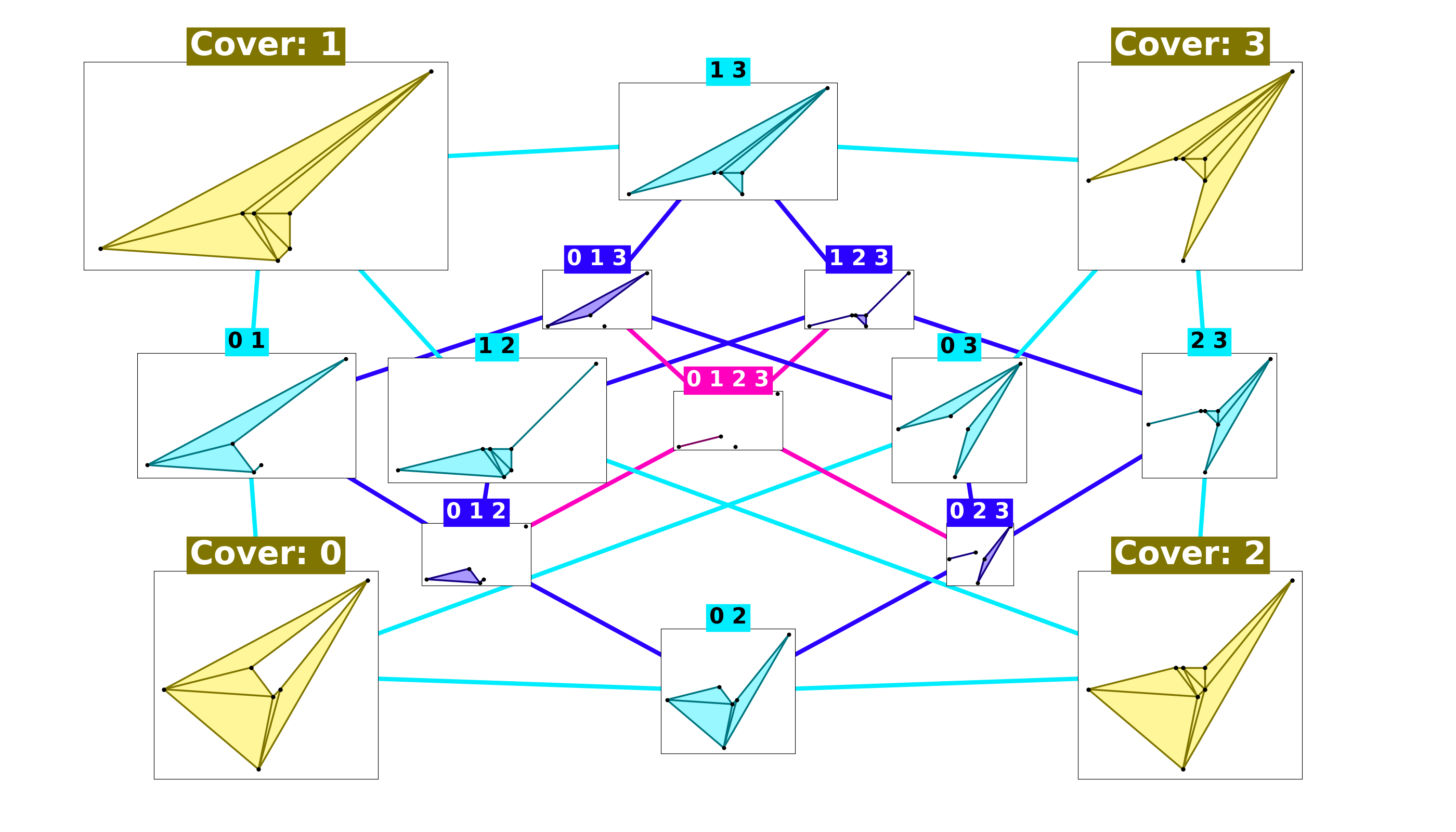}
    \caption{Covering sets over the nerve.}
    \label{subfig:cover-second-diff}
    \end{subfigure}
    \caption{Cover where there is a nontrivial second differential.}
    \label{fig:second-diff}
\end{figure}

\begin{figure}
    \centering \includegraphics[width=0.95\textwidth]{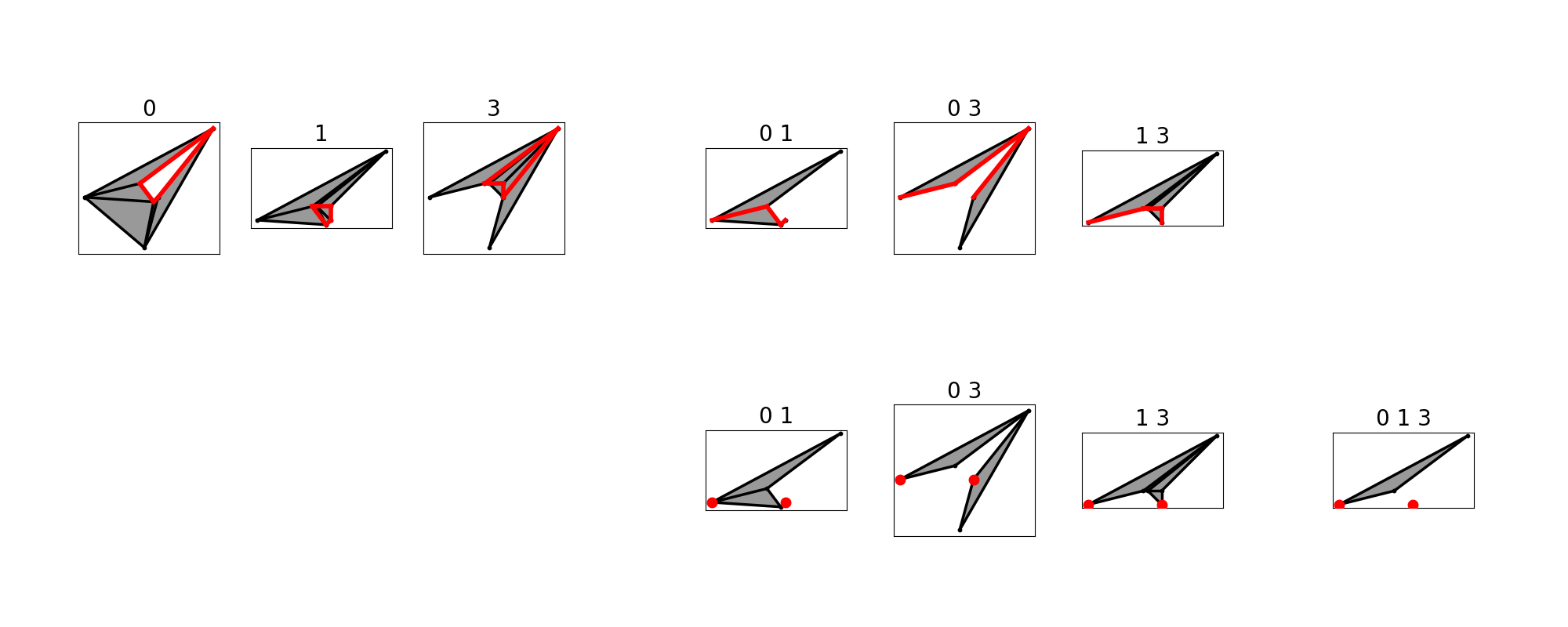}
    \caption{Illustration of a class representative involving the nontrivial second page differential $d_{2,0}^2$. We only depict the regions where such representative is nontrivial, i.e. the covers $0$, $1$ and $3$ and their intersections.}
    \label{fig:second-diff-reps}
\end{figure}

Computing the second page differentials adds algorithmic complexity which is unnecessary for the spectral sequences considered in this article. Unfortunately, there are some cases where such differentials are nontrivial, as shown next.
\begin{example}\label{ex:second-page-differential}
    Consider $8$ points divided into $4$ zones as shown in Subfigure~\ref{subfig:zones-second-page}. 
    This gives rise to the cover illustrated in Subfigure~\ref{subfig:cover-second-diff}. 
    Here, notice that the covering set $0$ has a hole, which means that there is an infinity interval within the first page term $E^1_{0,1}$. 
    Since this hole does not appear in any double intersection, there is no generator associated to an infinite interval in the domain of $d^1_{1,1}$ and so $E^2_{0,1}$ has a generator $\alpha$ associated to an infinite interval in its barcode.
    Such $\alpha$ is ``killed'' by the second page differential $d^2_{2,0}$. 
    To see this, consider the class represented by a couple of points in the covering set $[0,1,3]$ as depicted on the bottom right of Figure~\ref{fig:second-diff-reps}.
    Embedding this representative into the three covering sets $[0,1]$, $[0,3]$ and $[1,3]$, we see that these lift to three different paths shown in the middle of Figure~\ref{fig:second-diff-reps}. 
    Computing the image under the horizontal differential of these paths we obtain the three cycles shown on the left upper part of Figure~\ref{fig:second-diff-reps}. 
    Since one of these cycles is a representative for $\alpha$ while the other two are boundaries, this means that $\alpha$ lies in the image of $d^2_{p,q}$, and there is no infinite interval in the term $E^3_{0,1}$.
\end{example}
\leavevmode
\leavevmode

From Example~\ref{ex:second-page-differential}, one might deduce that infinity intervals play an important role in determining whether second page differentials vanish or not. 
In fact, we make the following claim which we will prove in a forthcoming article in more generality. 
\begin{claim}~\label{clm:collapse-second}
Suppose that  $E^2_{0,1}$, $E^2_{1,1}$ and $\Ker(d^1_{2,1})$ have no infinite intervals on their interval decompositions. Then, $E^2_{p,q}$ collapses on the second page.
\end{claim}

A practical consequence of Claim~\ref{clm:collapse-second} is that we only need to compute the persistent homologies up to triple intersections. 
If the hypotheses from~\ref{clm:collapse-second} holds, then we are ready to obtain the persistent homology of the alpha complex after computing the extension problem. 
In our experiments, we confirm Claim~\ref{clm:collapse-second} by checking that the resulting barcode obtained using the Mayer-Vietoris spectral sequence is equal (up to a tolerance of $10^{-8}$) to that obtained by following the sequential alpha complex persistence computation pipeline.
Further, Example~\ref{ex:second-page-differential} is an instance of a case where the hypotheses from Claim~\ref{clm:collapse-second} does not hold and the spectral sequence does not collapse on the second page.

\subsection{Extension Problem}\label{sub:extension-problem}
\begin{figure}
    \centering
    \begin{subfigure}[b]{0.3\textwidth}
    \centering
    \includegraphics[width=\textwidth]{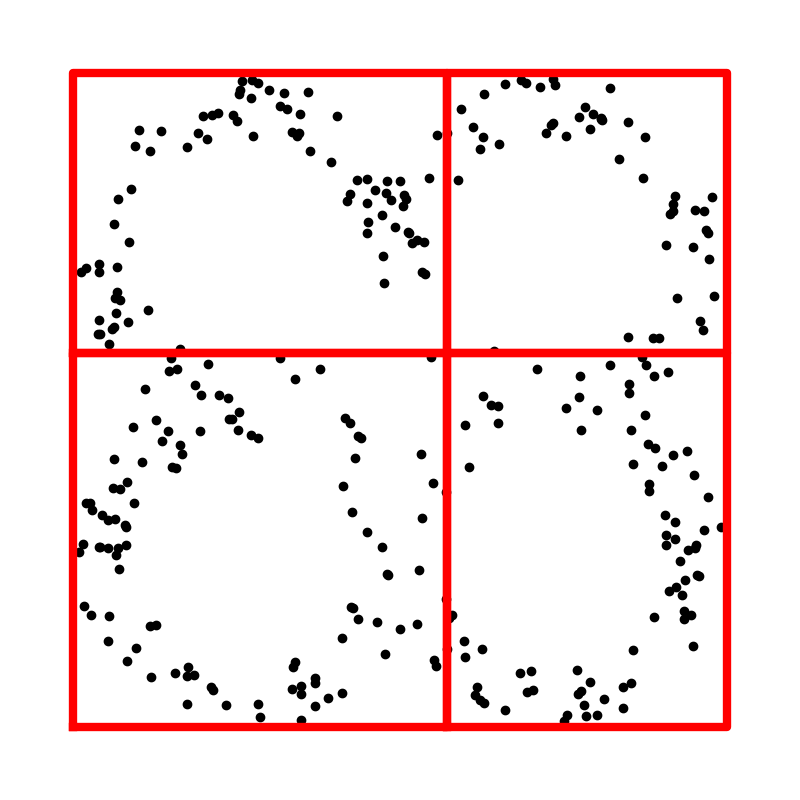}
    \caption{Cover zones.}
    \label{subfig:cover_extension}
    \end{subfigure}
    \hfill
    \begin{subfigure}[b]{0.33\textwidth}
        \centering
        \includegraphics[width=\textwidth]{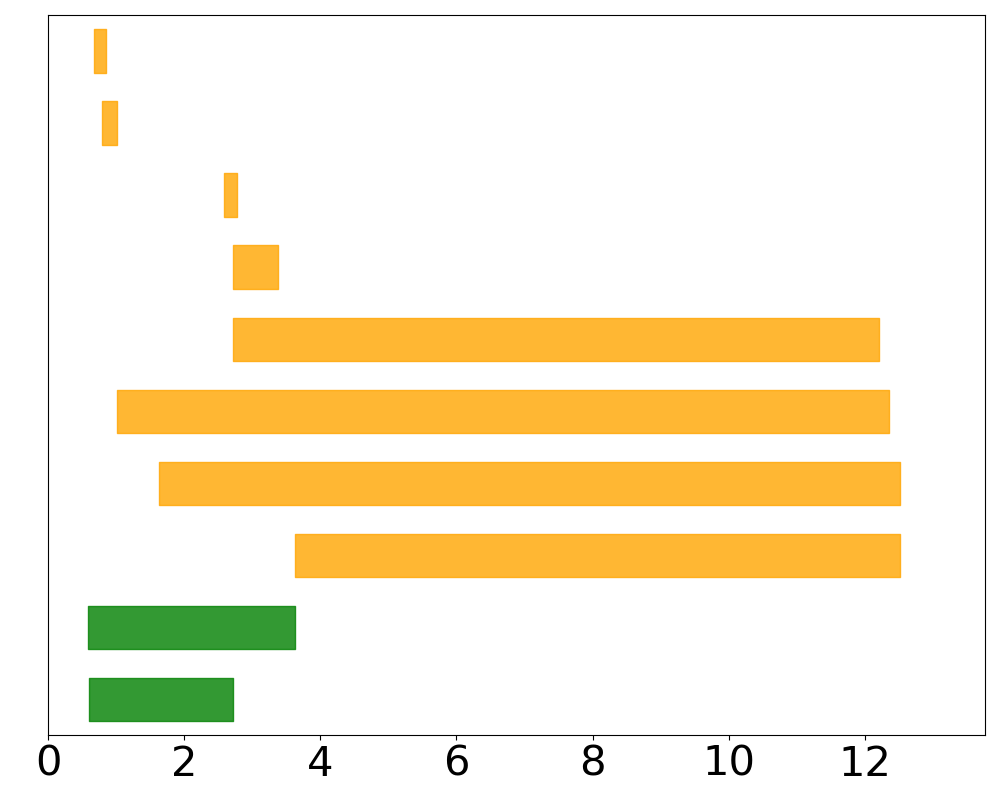}
        \caption{Barcodes on $E^2$ terms.}
        \label{subfig:E2-barcode-extension}
    \end{subfigure}
    \hfill
    \begin{subfigure}[b]{0.33\textwidth}
        \centering
        \includegraphics[width=\textwidth]{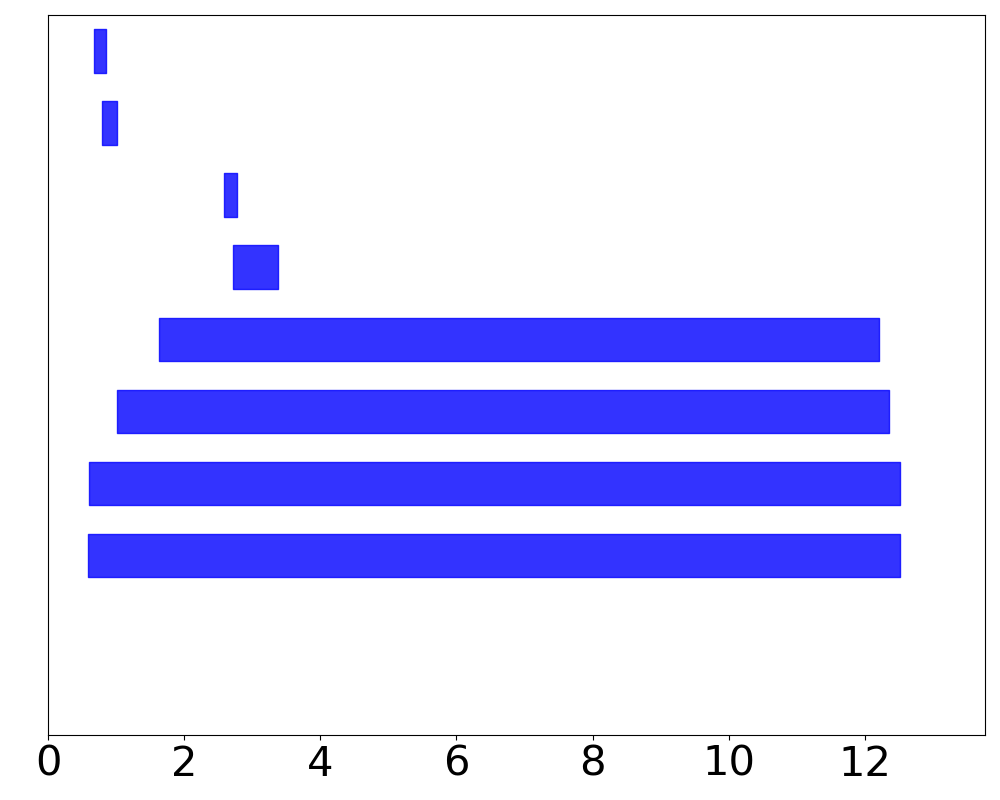}
        \caption{The $\PH_1$ barcode.}
        \label{subfig:PH1-barcode-extension}
    \end{subfigure}
    \caption{Zones and barcodes of an example where the extension problem is nontrivial. See Example~\ref{ex:extension_problem} for details.}
    \label{fig:extension_zones_barcode}
\end{figure}

In this section, we briefly describe how to solve the extension problem in our particular case by adapting the general setting from~\cite{Torras2023}. Suppose we have computed the second page of the spectral sequence $E^2_{p,q}$. Then, we can recover the $0$-dimensional persistent homology directly via the isomorphism $\PH_0(\AXn)\simeq E^2_{0,0}$. Recall that we assume that the spectral sequence collapses on the second page. Also, consider the total complex $C_m^\Tot=\bigoplus_{p+q=m} C_{p,q}$ on the double complex introduced in Subsection~\ref{sub:particular_permaviss} and recall that there is an isomorphism $\PH_1(\AXn)\simeq \PH_1(C^\Tot)$. In such case, there is a short exact sequence 
\begin{equation}\label{seq:ses_extension_alpha}
\begin{tikzcd}
0 \ar[r] &
E^2_{0,1} \ar[r] & 
\PH_1(C^\Tot) \ar[r] & 
E^2_{1,0} \ar[r] & 
0\,.
\end{tikzcd}
\end{equation}
Basically, the extension problem amounts to recovering $\PH_1(C^\Tot)$ from the terms $E^2_{0,1}$ and $E^2_{1,0}$ using the sequence~(\ref{seq:ses_extension_alpha}). This particular extension problem (using second page terms) was first solved in~\cite{YoonGhrist2020}. 
\begin{example}\label{ex:extension_problem}
    Consider a sample about four circles and the associated Mayer-Vietoris spectral sequence taking four zones, as depicted in Subfigure~\ref{subfig:cover_extension}. We compute the second page, obtaining the barcodes of the terms $E^2_{0,1}$ and $E^2_{1,0}$ which are depicted in orange and green respectively in Subfigure~\ref{subfig:E2-barcode-extension}. Here, we remark that the original barcodes had more intervals, but we are focusing on the ones that appear on the \emph{optimised entries} speedup (see Section~\ref{sec:implementation}) as well as those that are long enough. In this case, solving the extension problem amounts to recovering the intervals of $\PH_1(\AXn)$ depicted in Subfigure~\ref{subfig:PH1-barcode-extension} from the barcodes shown in Subfigure~\ref{subfig:E2-barcode-extension}.  
\end{example}
\leavevmode

In order to solve the extension problem, we need to take a closer look into the structure of the classes in $\PH_1(C^\Tot)$. 
A class $\alpha \in \PH_1(C^\Tot)$ with $\alpha\sim [a,b)$ is represented by a chain $(w_0, w_1) \in C_{0,1}\oplus C_{1,0}$ which is a cycle by the total differential $d^\Tot$. 
This means that\footnote{
Technically $d^V_1(w_0) \boxplus d^H_1(w_1) = Z_a$ as detailed in~\cite{Torras2023}.  For ease, we just write $0$ instead of $Z_a$.
} 
$d^V_1(w_0) \boxplus d^H_1(w_1) = 0$.
All these classes are given modulo boundaries by $d^\Tot_2$, that is, up to images $d^\Tot_2((a_0,a_1,a_2))$ for chains $(a_0,a_1,a_2) \in C_{0,2}\oplus C_{1,1}\oplus C_{2,0}$. 
From this viewpoint, $E^2_{0,1}$ is the submodule of $\PH_1(C^\Tot)$ whose classes can be represented by cycles $(w_0,w_1) \in C_{0,1}\oplus C_{1,0}$ with $w_1=0$. 
Sequence~(\ref{seq:ses_extension_alpha}) tells us that the term $E^2_{1,0}$ is isomorphic to $\PH_1(C^\Tot)$ modulo $E^2_{0,1}$. 
More concretely, a nontrivial class $\alpha \in E^2_{1,0}$ can only be represented by chains $(w_0,w_1) \in C_{0,1}\oplus C_{1,0}$ where $w_1\neq 0$ holds. 
At value $b$ the class $\alpha$ dies, which means that there exists a chain $(a_0,a_1,a_2) \in C^\Tot_2$ whose birth value is $b$ and such that $w_1 \boxplus d^V_1(a_1)\boxplus d^H_1(a_2)=0$. 
Here we point out that the barcode of $E^1_{2,0}$ only has intervals starting at $0$, so that $w_1$ and $w_1 \boxplus d_1^H(a_2)$ have the same birth value. In the implementation from Section~\ref{sec:implementation} we choose $w_1$ such that $a_2=0$.
\leavevmode
\begin{figure}
    \centering
    \includegraphics[width=0.9\textwidth]{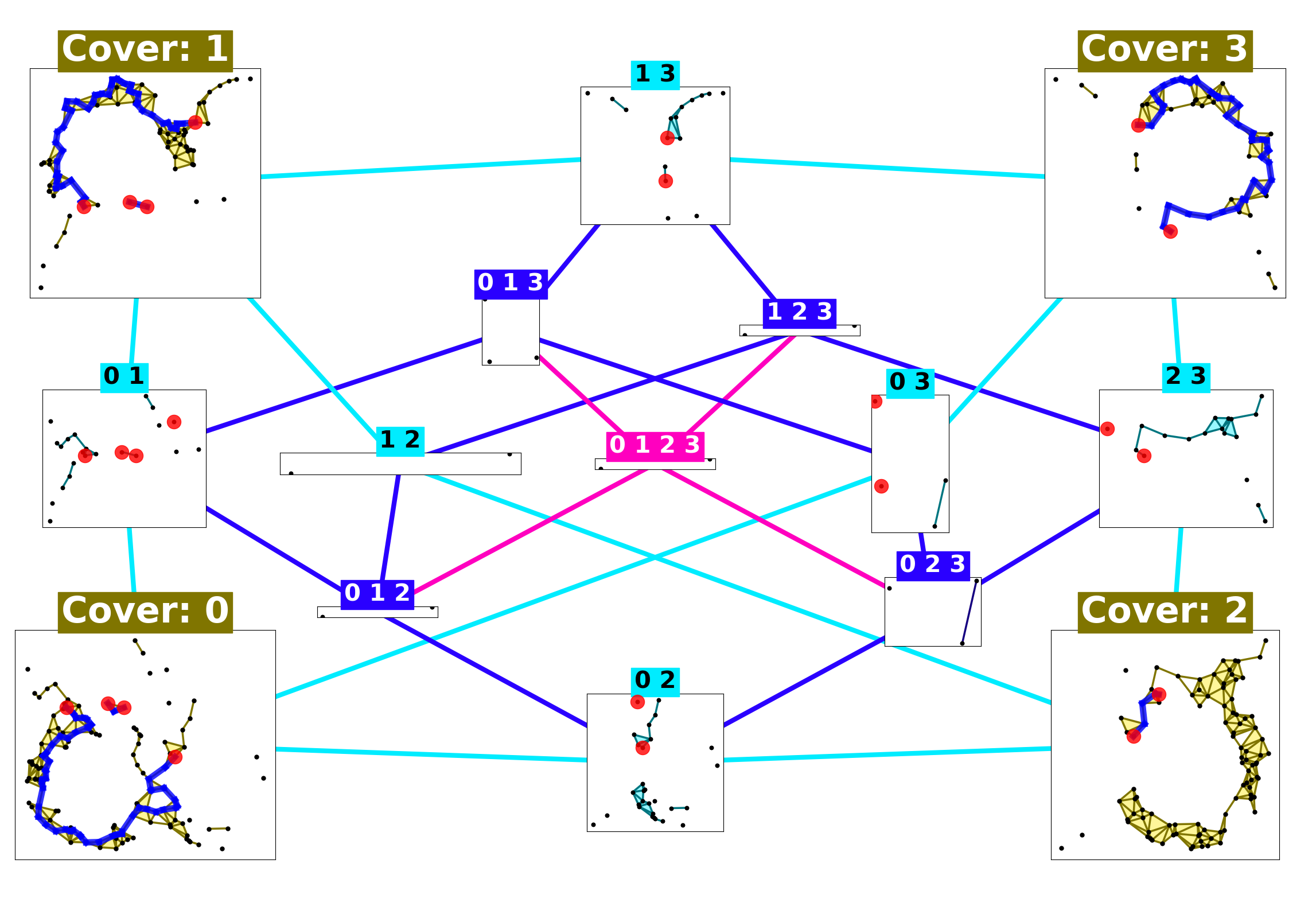}
    \caption{Representatives related to the birth of a class in $E^2_{1,0}$.}
    \label{fig:E2-1-0-class-birth}
\end{figure}

\begin{figure}
    \centering
    \includegraphics[width=0.9\textwidth]{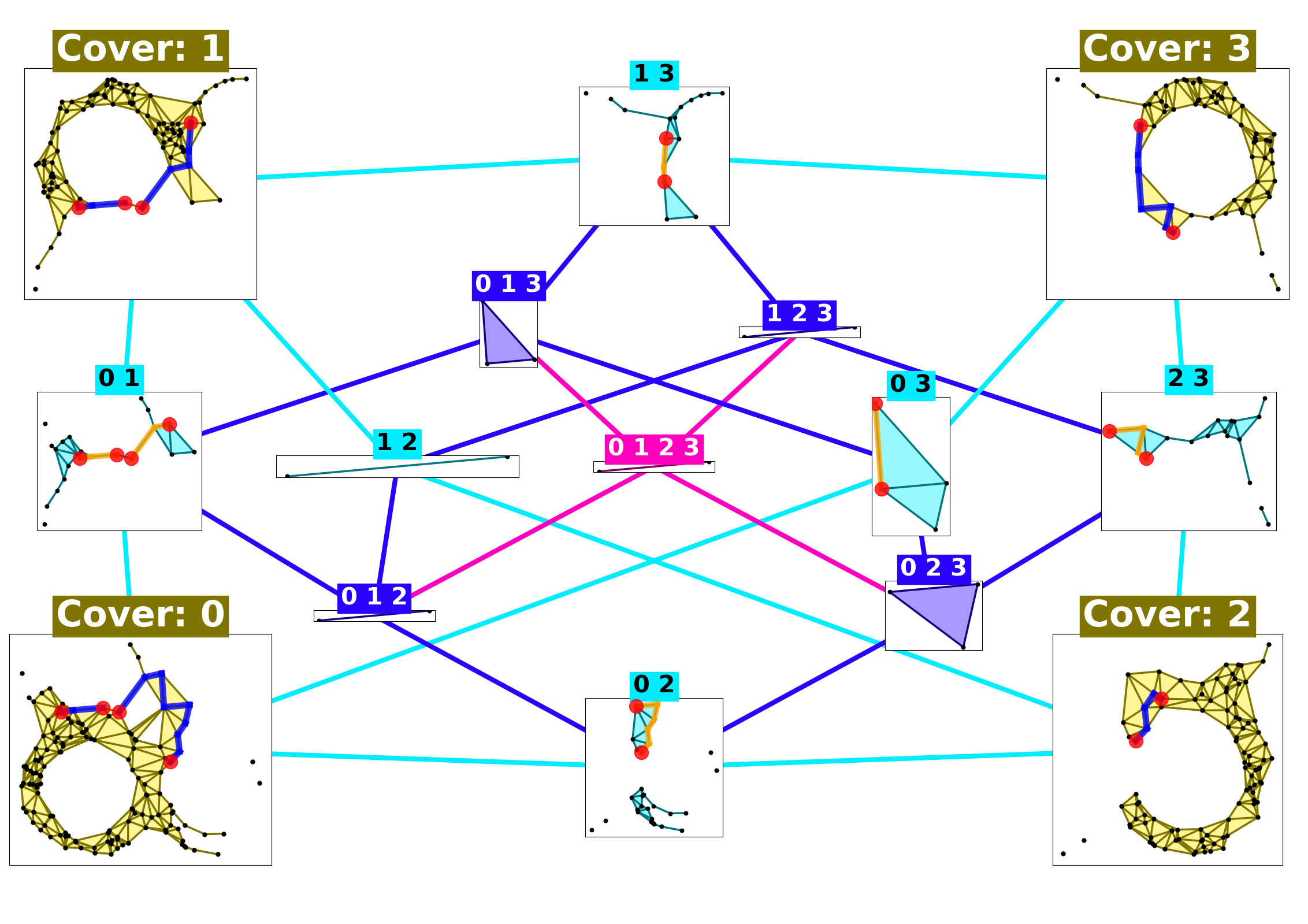}
    \caption{Representatives related to the death of a class in $E^2_{1,0}$.}
    \label{fig:E2-1-0-class-death}
\end{figure}

\begin{figure}
    \centering
    \begin{minipage}{.48\textwidth}
        \begin{subfigure}[b]{\textwidth}
            \centering
            \includegraphics[width=\textwidth]{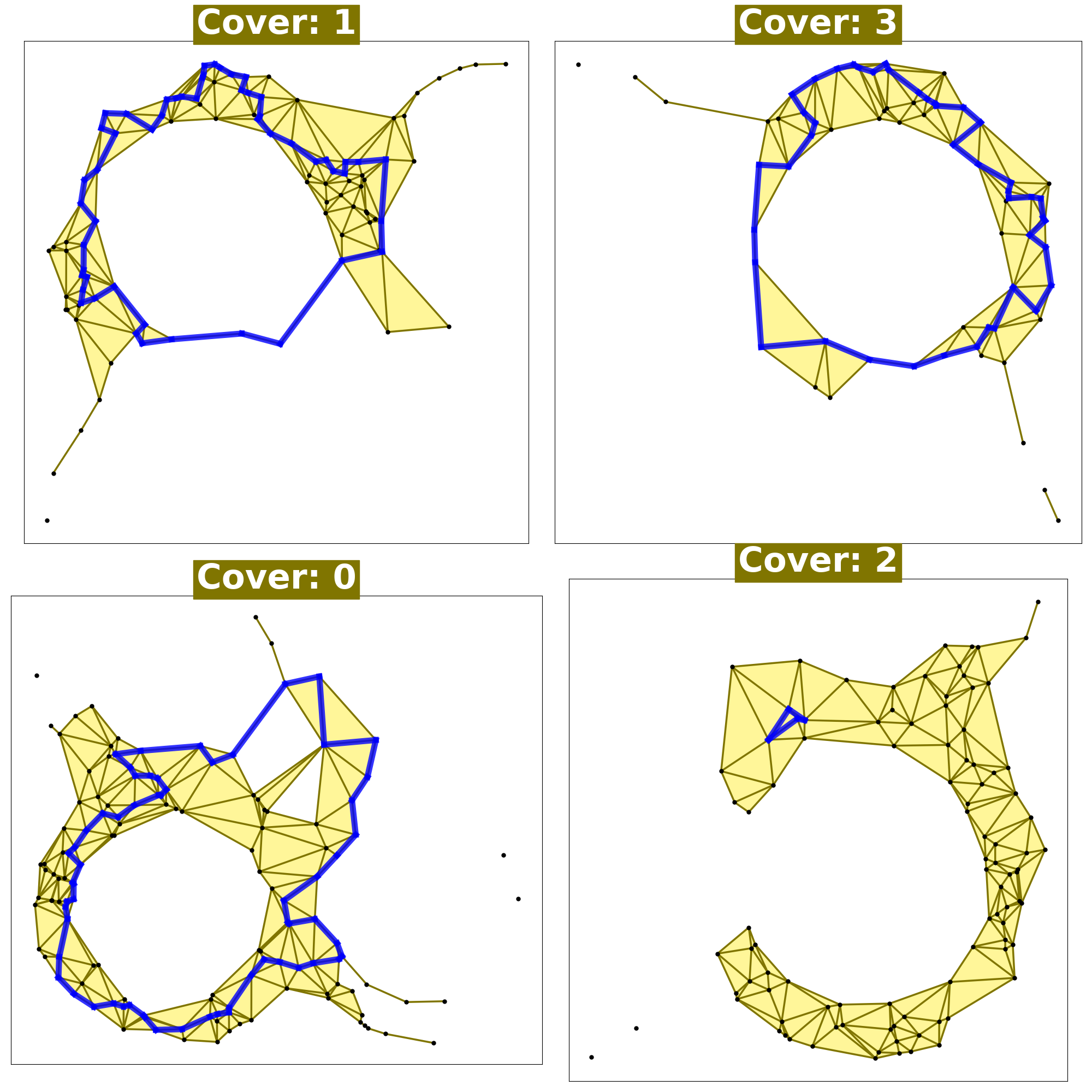}
            \caption{The chain $w_0\boxplus d_1^H(a_1)$.}
            \label{subfig:ext-reps}
        \end{subfigure}
    \end{minipage}
    \hfill
    \begin{minipage}{.48\textwidth}
    \begin{subfigure}[b]{\textwidth}
        \centering
        \includegraphics[width=\textwidth]{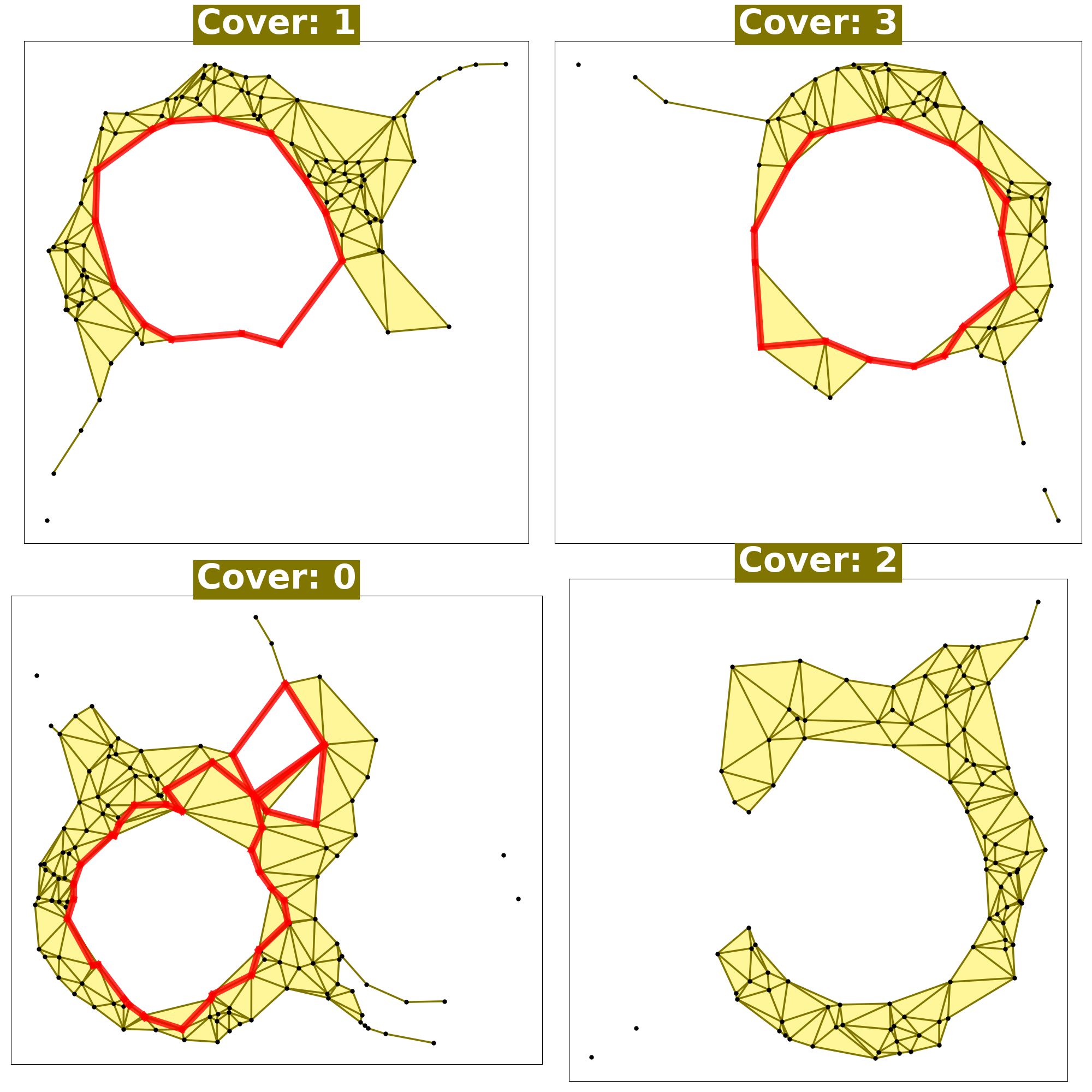}
        \caption{Local cycles}
        \label{subfig:local-cycles}
    \end{subfigure}
    \end{minipage}
    \caption{Extended cycles and their corresponding local cycles.}\label{fig:ext-local-cycles}
\end{figure}

\begin{example}~\label{ex:extension_reps}
    We continue Example~\ref{ex:extension_problem}.
    Consider the class $\beta_2\in E^2_{1,0}$ associated to the bottom green interval from Subfigure~\ref{subfig:E2-barcode-extension}; which we denote as $[a,b)$. This $\beta_2$ is represented by a chain $(w_0,w_1) \in C^\Tot_1$ where $w_1$ is given by the red dots from the double intersections of Figure~\ref{fig:E2-1-0-class-birth}. Although these red dots cannot be lifted on double intersections, the image $d^H_1(w_1) \in C_{0,0}$ given by the red dots on each cover (see Figure~\ref{fig:E2-1-0-class-birth} again) becomes a boundary of the paths indicated in blue, which lead to $w_0$. Now, since $\beta_2$ dies at $b$, there is a chain $(a_0,a_1,0) \in C^\Tot_2$ with birth value $b$ and such that $w_1 \boxplus d^V_1(a_1)=0$. We depict $a_1$ as the orange paths on the double intersections from Figure~\ref{fig:E2-1-0-class-death}. Next, we compute $d^H_1(a_1)$ which is depicted as the blue paths on the covers from Figure~\ref{fig:E2-1-0-class-death}. Thus, we have the equality in $\PH_1(C^\Tot)$
    \[
    \bOne_b(\alpha) = \bOne_b([(w_0, w_1)]) =\big[(w_0, w_1)  \boxplus d^\Tot_1(a_0,a_1,0)\big] = [(w_0 \boxplus d^H_1(a_1), 0)]\,.
    \]
    We can visualise the chain $w_0 \boxplus d^H_1(a_1)\in C_{0,1}$ as the blue cycles depicted on Subfigure~\ref{subfig:ext-reps}. Intuitively, these cycles are the result of ``gluing'' along the red dots the blue paths from Figure~\ref{fig:E2-1-0-class-birth} with the blue paths from Figure~\ref{fig:E2-1-0-class-death}.
\end{example}
\leavevmode
\leavevmode

We continue considering the class $\alpha\sim[a,b)$ in $E^2_{1,0}$. 
At filtration value $b$ the chain $(w_0\boxplus d_1^H(a_1), 0)$ represents a class in $E^2_{0,1}$. 
As we have already computed a barcode basis $\cE_0$ for the term $E^2_{0,1}$, we can write the class of $(w_0\boxplus d_1^H(a_1), 0)$ in terms of $\cE_0$.
Following the same procedure as in Appendix~\ref{app:associated-matrix} we obtain the \emph{extension coordinates} of $\alpha$ in terms of $\cE_0$.
Repeating this procedure for all generators of a basis $\cE_1$ of $E^2_{1,0}$, we compute the \emph{extension matrix} $\EXT_1(\cE_0, \cE_1)$ which is a $\vert \cE_0\vert \times \vert \cE_1\vert$ matrix whose columns contain the extension coordinates of the generators from $\cE_1$ in terms of $\cE_0$.

Subsequently, we can get a barcode basis for $\PH_1(C^\Tot)$ by taking a simple quotient.
For this, we consider the persistence module $\rB E^2_{1,0}$ that results from $E^2_{1,0}$ by turning all intervals into infinity intervals starting at the birth values. 
$\rB E^2_{1,0}$ has a barcode basis $\rB\cE_1$ which we can construct by taking generators $\alpha\sim [a,b)$ from $\cE_1$ and assigning them to infinity generators $\alpha^\rB\sim [a,\infty)$.
Similarly, we can construct the persistence module $\rD E^2_{1,0}$ by turning all intervals from $E^2_{1,0}$ into infinity intervals starting at the death values. 
The barcode basis $\rD\cE_1$ of $\rD E^2_{1,0}$ thus consists of generators  $\alpha^\rD\sim [b,\infty)$ which where constructed from generators $\alpha\sim [a,b)$ of $\cE_1$.

We can now consider the persistence morphism 
$\textrm{ext}_1: \rD E^2_{1,0} \longrightarrow E^2_{0,1}\oplus \rB E^2_{1,0}$ whose associated matrix in the bases $\rD \cE_1$ and $\cE_0\cup \rB\cE_1$ is given by the block matrix (where we assume that $\cE_1$, $\rB\cE_1$ and $\rD\cE_1$ all follow the same order)
\[
\left(\begin{array}{c} 
\EXT_1(\cE_0, \cE_1) \\
\hline
\Id
\end{array}\right)\,
\mbox{ which we use to compute }
\dfrac{E^2_{0,1}\oplus \rB E^2_{1,0}}{\textrm{ext}_1(\rD E^2_{1,0})}\,.
\]
This last quotient is isomorphic to $\PH_1(C^\Tot)$. To see this, consider a generator $\alpha \sim [a,b)$ from $\cE_1$ and notice that  $\textrm{ext}_1(\rD E^2_{1,0})$ identifies $\alpha^\rB$ with the extension coordinates of $\alpha$ at value $b$. Altogether, we obtain the barcode of $\PH_1(\AXn)$.
\begin{example}\label{ex:extension_matrix}
    We continue from Example~\ref{ex:extension_reps}. In particular, we consider the chain $w_0 \boxplus d^H_1(a_1)\in C_{0,1}$ which we obtained from a representative of $\beta_2$ and which is depicted as blue cycles on Subfigure~\ref{subfig:ext-reps}. We can write these cycles in terms of the basis on $E^2_{0,1}$. In particular, we obtain that the extension coordinates of $\beta_2$ are $(0,0,1,1,1,1,1,0)$; where the nontrivial entries lead to the red cycles from Subfigure~\ref{subfig:local-cycles}. We repeat this procedure for the class $\beta_1$ from $E^2_{1,0}$, which is associated with the top green interval from Subfigure~\ref{subfig:E2-barcode-extension}. In this case, the extension coordinates are  $(0,0,0,0,0,0,0,1)$. We can now obtain the barcode from~\ref{subfig:PH1-barcode-extension} by computing a quotient. For completeness, we compute this quotient in Appendix~\ref{app:example-quotient}.
\end{example}

\section{Implementation Details}\label{sec:implementation}

We implemented the algorithm described 
above in C++ on top of Open MPI~\cite{mpi}.
In this section, we give an overview of the libraries we used as well as some technical details.

\subsection{Program Structure}\label{sub:program-structure}
\begin{figure}
\centering
\begin{tikzpicture}
    \node[draw, circle, line width=0.05cm, inner sep=1pt] (R0) at (-1,2) {$D_0$};
    \node[draw] (L0) at (1,2) {$AC_0$};
    \node[draw] (L1) at (1,1) {$AC_1$};
    \node[draw] (L2) at (1,0) {$\cdots$};
    \node[draw] (LM) at (1,-1) {$AC_{M-1}$};
    \node[draw] (P0) at (2.75,2) {$PH_0$};
    \node[draw] (P1) at (2.75,1) {$PH_1$};
    \node[draw] (P2) at (2.75,0) {$\cdots$};
    \node[draw] (PM) at (2.75,-1) {$PH_{M-1}$};
    \node[draw] (S0) at (5,2) {$E2_0$};
    \node[draw] (S1) at (5,1) {$E2_1$};
    \node[draw] (Lf0) at (7,2) {$L_0$};
    \node[draw] (Lf1) at (7,1) {$L_1$};
    \node[draw] (Lf2) at (7,0) {$\cdots$};
    \node[draw] (LfM) at (7,-1) {$L_{M-1}$};
    \node[draw] (Ex1) at (9,1) {$Ex_1$};
    \node[draw, circle, line width=0.05cm, inner sep=1pt] (G0) at (9,2) {$G_0$};
    \draw (R0.east)--(L0.west);
    \draw (R0.east)--(L1.west);
    \draw (R0.east)--(L2.west);
    \draw (R0.east)--(LM.west);
    \draw[Stealth-Stealth] (L0.south) -- (L1.north);
    \draw[Stealth-Stealth] (L0.south) to [out=350,in=10] (L2.north);
    \draw[Stealth-Stealth] (L1.south) -- (L2.north);
    \draw[Stealth-Stealth] (L1.south) to [out=340,in=20] (LM.north);
    \draw[Stealth-Stealth] (L2.south) -- (LM.north);
    \draw[Stealth-Stealth] (L0.south) to [out=200,in=160] (LM.north);
    \foreach \i in {0,1,2,M} {
            \draw (L\i.east)--(P\i.west);
    }
    \foreach \i in {0,1,2,M} {
        \foreach \j in {0,1} {
            \draw (P\i.east)--(S\j.west);
        }
    }
    \draw (S1.east)--(Lf1.west);
    \draw[dashed] (P2.east)--(Lf2.west);
    \draw[dashed] (PM.east)--(LfM.west);
    \foreach \j in {0,1,2,M} {
        \draw (S0.east)--(Lf\j.west);
    }
    \draw (Lf0.east)-- (Ex1.west);
    \draw (Lf1.east)-- (Ex1.west);
    \draw (Lf2.east)-- (Ex1.west);
    \draw (LfM.east)-- (Ex1.west);
    \draw[-stealth] (Ex1.north)--(G0.south);
    \draw (Lf0.east)--(G0.west);
    \draw (Lf2.east)--(G0.west);
    \draw (LfM.east)--(G0.west);
    
\end{tikzpicture}
\caption{Computation workflow for PerMaViss++.}
\label{fig:permaviss_workflow}
\end{figure}
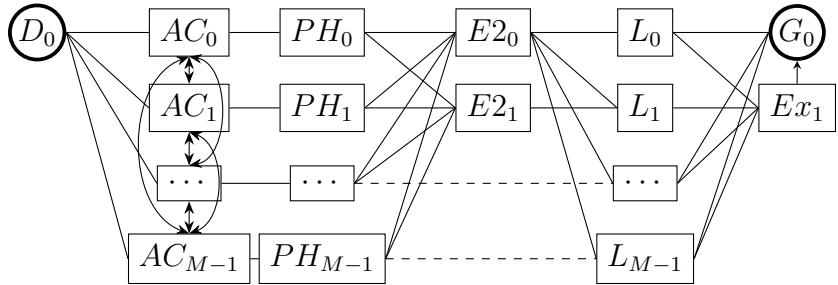

The program is structured as shown in Figure~\ref{fig:permaviss_workflow}. Recall that $M$ denotes the number of processes. The program runs the following steps: 
\begin{enumerate}
    \item ($D_0$) First, processor $0$ takes the input data, computes the bounding boxes, and, subsequently, distributes the points of the respective zone as well as the grid data to all processors.
    \item ($AC_i$) Next, each processor computes its alpha complex, doing the necessary information exchanges with other processes (additional points, critical edges/vertices and critical non-Gabriel edges)
    \item ($PH_i$) Next, processor $i$ computes the persistent homology of cover $\rA_i$ as well as intersections $\rA_i\cap \rA_j$ and $\rA_i\cap \rA_j \cap \rA_k$. In addition, processor $i$ also computes the matrices associated with inclusions from double and triple intersections. Barcodes and matrices are sent to processors $0$ and $1$. To reduce the size of these, the \emph{optimised entries speedup} is used.
    \item ($E2_i$) Process $i$ handles $i$-dimensional persistent homology barcodes and matrices, for $i=0,1$. 
    In this case, the first two rows of the second page $E^2_{p,q}$ are computed. After, process $0$ sends barcodes and representative information related to the term $E^2_{1,0}$ to other processes.
    \item ($L_i$) Process $i$ handles representative chains in $C_{1,0}$ and lifts them to chains in $C_{0,1}$, computing also their classes in $E^1_{0,1}$. This leads to blocks of the matrix $\EXT_1(\cE_0, \cE_1)$ related to each process.
    \item ($Ex_1$) Process $1$ handles the extension matrix and computes the barcode of $\PH_1(\AXn)$ restricted to the optimised entries. This barcode is sent to process $0$.
    \item ($G_0$) Process $0$ gathers all the intervals which were not used due to the optimised entries approach. Optionally, it also compares its resulting barcode with the result returned by other software, such as GUDHI~\cite{gudhi}.
\end{enumerate} 
\leavevmode

\subsection{Parallel Alpha Complex Computation with \texttt{Alpha\_p}}
The construction of the Delaunay triangulation and the computation of the Alpha complex filtration values as described in Section~\ref{sec:alpha-construction} is handled by the class \texttt{Alpha\_p} (steps ($D_0$) and ($AC_i$) in Subsection~\ref{sub:program-structure}).
It is based on the \texttt{Alpha\_complex} class from the GUDHI library~\cite{gudhi}.
The constructor of \texttt{Alpha\_p} creates the Delaunay triangulation of a two-dimensional point cloud given a file containing its coordinates of a two-dimensional point cloud.
As already mentioned in Subsection~\ref{sub:alpha-covers}, we follow the algorithm described in \cite{Lo2012} to construct the Delaunay triangulation in parallel.
The computation of the filtration values is then handled by the member function \texttt{create\_simplex\_trees}.

\subsubsection{\texttt{Alpha\_p} constructor}

When the constructor of \texttt{Alpha\_p} is called, the processor with rank 0 reads in the data from a file and computes the corresponding global bounding box $B$.
Next, processor 0 computes the decomposition of the bounding box into a set of smaller cells based on a rectilinear grid that depends on the density parameter which is set by the user.
This parameter describes the average number of points lying in a grid cell.
Cells are then grouped into bigger zones by processor 0 based on a coarser rectilinear grid that depends on the number $M$ of cores available.
Currently, our implementation requires that $M=m_1 \times m_2$ for some $m_1, m_2 \in \mathbb{N}$.
Subsequently, processor 0 distributes the data to the other cores, as well as the coordinates of the local bounding boxes $B_i$, the coordinates of the global bounding box $B$ and the cell size.
Eventually, processor 0 only stores its own sub-point cloud $\bX_0$.
This completes step ($D_0$).

After the data distribution, each core computes its local Delaunay triangulations $\DXi{i}$ using the implementation provided by the CGAL library~\cite{cgal}.
Iteratively, the triangulations are then expanded in parallel as described in Subsection \ref{sub:alpha-covers}.
In each iteration, data points are exchanged between the processors.
This is the first part of step ($AC_i$).

\subsubsection{The \texttt{create\_simplex\_trees} function}
To store the simplicial complexes $\rA_i$, $i=0,\ldots,M-1$, and their intersections, we use the \texttt{Simplex\_tree} class from the GUDHI library.
We use the function \texttt{create\_complex} from the alpha complex class as a model for our modified version \texttt{create\_simplex\_trees}.
The same holds true for the other helper functions needed in this step.

When the function is called, it first creates the simplex tree for the cover by going through all triangles in the Delaunay triangulation and only storing the ones that are either inner or boundary triangles.
For each boundary triangle, we compute the one or two other intersecting zones and we add the triangle to the simplex tree storing the corresponding two- and three-fold intersection.
As explained in Subsection \ref{sub:alpha-intersection}, the three-fold intersections are then shared with the other cores and, subsequently, the two- and three-fold intersections are updated with the critical vertices and edges.
Next, the function then computes the filtration values for $\rA_i$ while keeping track of all detected non-Gabriel edges.
In order to do this, we again modified the necessary helper functions from the GUDHI library.
Afterwards, the filtration values for the intersections are computed.
As described in Subsection~\ref{sec:alpha-fv}, we have to be extra careful with the non-Gabriel edges.
Thus, for each edge, we check whether it is in the list of non-Gabriel edges computed before.
If so, the filtration value is adopted. 
In those cases, we also run a check whether it would have been detected in the intersection as a non-Gabriel edge.
Those that were not detected are stored as critical non-Gabriel edges.
When the computation of the filtration values of the intersections is done, the lists of critical non-Gabriel edges are shared with the other processors and the filtration values of the non-Gabriel edges are updated.
This completes step ($AC_i$).

\subsection{Spectral Sequence Computatation}
In this subsection, we give some guidelines for the implementation of the Mayer-Vietoris spectral sequence. 
By this point, each process has computed the simplex trees associated to its cover element $\rA_i$ as well as double intersections $\rA_i\cap \rA_j$ and triple intersections $\rA_i \cap \rA_j \cap \rA_k$ ranging over $j,k\neq i$. In particular, process $i$ is now ready to execute step $(PH_i)$ from Subsection~\ref{sub:program-structure}.

\subsubsection{First page computation}

The computation of the first page terms $E^1_{p,q}$ is done in parallel. Each processor computes the persistent homology in dimensions $0$ and $1$ on the cover element, the double intersections and the triple intersections that it is handling. 
Since we need to have access to the representatives and, in later steps, the reduced boundary matrix, we perform this computation using the standard reduction algorithm from PHAT~\cite{phat}. 
Notice that currently PHAT only handles computations on $\bZ_2$, which for our case is enough as mentioned in Subsection~\ref{sub:PH-alpha}. 
At some stages in the algorithm we need to track the column reductions performed by PHAT, such as when computing the matrix associated to a persistence morphism or when computing the chain whose boundary is a given cycle. 
For this, we had to add this functionality to the current code from PHAT.

The next step, for processor $i$, is to create barcode basis objects for the computed persistent homology groups. 
Also, the process computes the inclusion matrices $D^1_{m,[i],[i,j]}$ associated to $\PH_m(\rA_i\cap \rA_j) \hookrightarrow \PH_m(\rA_i)$ as well as inclusions from triple intersections to double intersections. 
After this is finished, the computed interval decompositions and the induced matrices are sent to processors $0$ and $1$, respectively, which handle the $0$ and $1$ rows from the first page of the spectral sequence.
Here we use the optimised entries approach, which we explain now. 
This optimisation is used since the terms $\PH_0(\rA_i)$ and $\PH_1(\rA_i)$ are too large and sending all their generators hinders the efficiency of the distribution of computations. 
Now, processor $i$ does not send the generators from $\PH_m(\rA_i)$ such that their associated rows in all matrices $D^1_{m,[i],[i,j]}$ are zero.
The intervals associated with such generators are not shared until the last step ($G_0$) when process $0$ claims all this withheld information to obtain the interval decomposition from $\PH_m(\AXn)$.
The reason why this optimisation works is that these reserved entries do not play any role when computing images or kernels of the differentials on the first page. 
Some of these entries might need to be shared during the computation of the extension problem, which we will comment in Subsection~\ref{subb:extension_problem_implementation}.
This completes step ($PH_i$).

\subsubsection{Second page computation with \texttt{E2\_page}}

Here, we describe the constructor for the class \texttt{E2\_page} which computes the terms of the $E^2$-page.
Processors $0$ and $1$ receive the information from other processors to obtain the barcodes from the first page terms, as well as the associated matrices.
When the data is received by the processors 0 and 1, respectively, the matrices corresponding to the first page differentials are filled.
Note that the matrices are stored as sparse matrices.
Since we compute persistent homology with coefficients in $\bZ_2$,
for each column, we only store the indices of the non-zero elements.
Next, columns and rows are reordered following the standard and endpoint order, respectively.
The resulting matrices are reduced using the function \texttt{image\_kernel} as described in~\cite{Torras2023} to compute the images and kernels of the maps. Afterwards, some quotients are computed by using the box gauss reduce algorithm from~\cite{Torras2023}; this finishes steps $(E2_i)$\fr{.} 

At this stage, processor $1$ checks that the hypotheses from Claim~\ref{clm:collapse-second} holds so that the spectral sequence collapses by the second page. If the check fails, it raises an error and the program stops. 
To prepare for the extension problem, core 0 sends the barcode and first-page representative coordinates from the $E_{1,0}^2$ term to other processors. 

\subsubsection{The Extension Problem}\label{subb:extension_problem_implementation}

The extension problem is handled in two steps. The first one is part $(L_i)$
and is done in parallel. 
More concretely, processor $i$ receives from processor $0$ the interval decomposition of the term $E^2_{1,0}$ as well as coordinates of representatives restricted to the summands $\PH_0(\rA_i\cap \rA_j)$ for a varying $j\neq i$. 
Now, for each generator of $E^2_{1,0}$, processor $i$ proceeds similarly as in Example~\ref{ex:extension_reps} and computes a cycle representative in $C_1(\rA_i)$. 
For each such representative, processor $i$ computes its coordinates in terms of the already stored barcode basis for $\PH_1(\rA_i)$. 
All the obtained coordinates are then sent to processor $1$. 
Recalling the optimised entries speedup, it is worth mentioning that there might be generators from the basis of $\PH_1(\rA_i)$ that were not previously shared with core $1$ but now come up as nontrivial coordinates. 
If this is the case, these generators are also shared with processor $1$. 

The second step to solve the extension problem ($Ex_1$) is completely handled by processor $1$.
First, it receives all coordinates from extended cycles so that it obtains representatives in the term $E^1_{0,1}$ for all generators from $E^2_{1,0}$. 
Next, processor $1$ computes the classes of these representatives in the term $E^2_{0,1}$, so that it obtains the matrix $\EXT_1(\cE_0, \cE_1)$.
With all this information, processor $1$ can proceed and compute $\PH_1(\AXn)$ as explained in Subsection~\ref{sub:extension-problem}.

\subsubsection{Gather and compare}
The persistent homology is now completely computed but stored and managed by different cores.
For a clear output, the function \texttt{gather} from the \texttt{E2\_page} class collects all persistent pairs into one sorted list per dimension.
This is step ($G_0$) in Subsection~\ref{sub:program-structure}.
When \texttt{gather} is called, all persistent pairs that were not in the images of the first page differentials,
and thus are not necessarily stored by core 0 or 1, are sent to core 0.
Core 1 also sends the persistent pairs resulting from the solution of the extension problem to core 0.
Core 0 then adds those to the barcode bases for dimension 0 and 1, respectively, and, eventually, sorts them according to the standard order.

If required, core 0 also computes persistent homology using GUDHI and compares the barcodes obtained by both methods.

\subsection{Experiments }

We have tested our implementation with two different artificial datasets consisting of $10^6$ data points. 
In the first sample the points are uniformly distributed on a rectangle and in the second one they are distributed on a circle with some noise.
Both complexes consist of about $6.0\times 10^6$ simplices in total.
We have tested our implementation on a NUMA machine with 48 Intel Xeon E7-8870 CPUs.
For the experiments, we chose first a very small density of 30 points per cell and then a density of 1000 points per cell. 
For each experiment, we computed the mean runtime of three runs. 
The results are shown in Table~\ref{tab:runtimes}.
\leavevmode
\begin{table}[]
    \centering
    \begin{tabular}{||c|c|c|c|c|c|c|c|c||}
        \hline
        \#processors & & 1 & 2 & 4 & 8 & 9 & 16 & 32 \\ [0.3ex]
        \hline\hline
        circle ($\rho = 30$) & total time &  - & 1869s & 777s & 488s & 535s & 430s & - \\
        circle ($\rho = 1000$) & total time & 396s & 1416s & 681s & 579s & 440s & 380s & 687s \\
        \hline
         & Delaunay & 91s & 76s & 65s & 69s & 55s & 51s & 69s \\
         & $\alpha$-complex & 192s & 918s & 464s & 226s & 220s & 164s & 149s \\
         & E1 & 203s & 98s & 48s & 23s & 25s & 16s & 9s \\
         & E1-diff & - & 11s & 18s & 22s & 20s & 27s & 19s \\
         & E2 & - & 16s & 19s & 41s & 51s & 119s & 279s \\
         & extension & - & 359s & 129s & 283s & 150s & 163s & 433s \\
         \hline
        uniform ($\rho = 30$) & total time & - & 3556s & 1964s & 710s & 605s & 521s & 965s \\
        uniform ($\rho = 1000$) & total time & 397s & 2107s & 1007s & 694s & 638s & 463s & 1008s \\
        \hline
        & Delaunay & 88s & 74s & 65s & 279s & 41s & 30s & 41s \\
         & $\alpha$-complex & 187s & 1459s & 720s & 279s & 277s & 126s & 83s \\
         & E1 & 208s & 101s & 48s & 21s & 20s & 10s & 5s \\
         & E1-diff & - & 11s & 18s & 14s & 20s & 12s & 6s \\
         & E2 & - & 19s & 26s & 72s & 80s & 150s & 472s \\
         & extension & - & 509s & 196s & 346s & 285s & 289s & 765s \\
        \hline        
    \end{tabular}
    \caption{Permaviss++ runtimes ($\rho$ denoting the chosen density)}
    \label{tab:runtimes}
\end{table}
\leavevmode
\begin{table}[]
    \centering
    \begin{tabular}{||c|c|c|c|c||}
        \hline
        program & sample & \#processors & mean (3 runs) \\ [0.3ex]
        \hline\hline
        GUDHI &  uniform & 1 & 199s \\
         & circle & 1 & 199s \\
         \hline
        PHAT + & uniform & 1 & 378s \\
        GUDHI alpha complex & circle & 1 & 373s \\
        \hline
        Permaviss++ & uniform & 1 & 397s \\
        (density = 1000) & circle & 1 & 396s \\
         & uniform & 16 & 463s \\
         & circle & 16 & 380s \\ [1ex]
        \hline        
    \end{tabular}
    \caption{Comparison of runtimes with GUDHI and PHAT}
    \label{tab:comparison_gudhi_phat}
\end{table}

We have further tested what slows down Permaviss++.
The runtime is currently mainly taken by three functions: 
The solution of the extension problem, the E2-page computation and the computation of the filtration values.
For the latter we have already discussed ways to speed up the computation as, currently,
the computations for the cover and then for the intersections are done sequentially.
We plan to compute them in parallel in the next version.
Another improvement we are currently working on is further reducing the size of the intersections by not including the whole boundary triangles but only certain edges of those.
This would also reduce the dimension of the intersections to one.
During the extension problem, the computation of lifts of generators to extend, i.e. steps $(L_i)$ from Subsection~\ref{sub:program-structure} takes most of the time. 
Additionally, in some runs processor $1$ spends a lot of time waiting for other processors to finish those computations.
Since this is a computing step handled by many processors simultaneously, we hope to improve it in the future. 
Another part that we think might be slowing computations is the current implementation of the box gauss reduce function from~\cite{Torras2023} in the extension problem computation. 
For future improvements, one might use the fact that the biggest block in such a reduction is a permutation matrix, as illustrated in Appendix~\ref{app:example-quotient}.

\section{Conclusion}\label{sec13}

We have presented a new algorithm that distributes the computation of persistent homology for alpha complexes on points in $\bR^2$. 
Our approach has the added benefit that one has access to more information than the standard persistent homology computation. This information includes localised homology~\cite{ZomCar2008} with respect to the cover that we consider; i.e. which intervals come from each cover element as well as which intervals result from the combination of the covering regions. We would like to explore use cases of this additional information in the future; notice that there are stability bounds for such information~\cite{TorrasPennig21}.

We believe that there is much room to optimise the algorithm presented here so that its parallelisation is more effective. For example, 
we intend to reduce the size of double and triple intersections so that also their dimensions are reduced to 1. 
This would significantly decrease the computation time of the second page terms and the extension problem. 
Also, our code is in the early stages and we hope to optimise it in the future.

Our next research direction will be to work on the article to show Claim~\ref{clm:collapse-second}, which will be more theoretical in nature than this work. Also, in the near future we also plan to work on the spectral sequence associated to three-dimensional Alpha complexes. Since Lo's algorithm works for three dimensions~\cite{Lo2012}, we are confident we can apply the experience from this work to tackle this future case which is more challenging.

\section*{Acknowledgements}

F. Jensen's research has been supported by the Deutsche Forschungsgemeinschaft (DFG, German Research Foundation) under Germany’s Excellence Strategy EXC 2181/1 - 390900948 (the Heidelberg STRUCTURES Excellence Cluster), Exploratory Project 5.1. 
The authors thank the Heidelberg STRUCTURES Excellence Cluster for financial support.
\'A. Torras-Casas research has been supported by the Engineering and Physical Sciences Research Council (EPSRC) grant EP/W522405/1 and Ministerio de Ciencia e Innovación grant with reference TED2021-129438B-I00. 

The authors further would like to thank G. Kanschat for both numerous fruitful discussions along the project as well as financial support.
The authors would further like to thank S. Ospina de los Rios for his support with all kinds of questions revolving around the implementation.

\section*{Declarations}

\begin{itemize}
\item We will publish the associated code as soon as possible. 
\item Both F. Jensen and \'A. Torras-Casas have equally contributed to this work.
\end{itemize}

\noindent

\printbibliography

\appendix

\section{Computing associated matrices}\label{app:associated-matrix}

Consider an inclusion of filtered simplicial complexes $f:L_*\hookrightarrow K_*$. In particular, we are interested in the persistence morphisms $f^*:\PH_n(L_*)\rightarrow \PH_n(K_*)$ that one obtains for all $n \geq 0$. As mentioned in Subsection~\ref{sub:persistence-morphisms}, if we have the associated matrix $F$ for $f^*$, then we can use it to compute barcode bases for the image and kernel of $f^*$. Here we explain briefly how to obtain such a matrix. 

Here we assume that our method for computing persistent homology leads to cycle representatives. Now, when computing $\PH_n(L_*)$, we obtain a (non-unique) set of cycle representatives $\Cycle_n(L_*)=\{w_i^L \in S_n(L_{a_i}) \colon d_n(w_i^L)=0\}_{i=1}^N$ such that their respective homology classes determine a barcode basis $\cA=\{\alpha_i\sim [a_i,b_i)\}_{i=1}^N$ for $\PH_n(L_*)$; where $\alpha_{ia_i}(1_{\bZ_2}) = w_i^L + \Ima(d_{n+1})_{a_i}$ for all $i=1,\ldots, N$. In the same way, one can compute $\PH_n(K_*)$ and obtain a (non-unique) set of cycle representatives $\Cycle_n(K_*)=\{w_i^K \in S_n(K_{a_i}) \colon d_n(w_i^K)=0\}_{i=1}^M$ such that their respective homology classes determine a barcode basis $\cB=\{\beta_i\sim [c_i, d_i)\}_{i=1}^M$. Using the inclusion $L_*\hookrightarrow K_*$, one can embed cycles $w^L_i$ from $\Cycle_n(L_*)$ to obtain cycles $\widetilde{w}_i^L$ in $K_*$. 
Given $1 \leq i \leq N$, we define the set $R_i=\{1 \leq j \leq M: a_i \in [c_j, d_j)\}$, and consider the following equalities 
\begin{multline*}
f^*(\alpha_i)_{a_i}(1_{\bZ_2}) = f^*_{a_i} \alpha_{ia_i}(1_{\bZ_2}) = f^*_{a_i}\big(\widetilde{w}_i^L + \Ima(d_{n+1})\big) 
\\ \stackrel{(*)}{=} \sum_{j \in R_i} F_{i,j} w_j^K + \Ima(d_{n+1}) = \sum_{j \in R_i} F_{i,j}\beta_{ja_i}(1_{\bZ_2}) = \bOne_{a_i}\bigg( \bigBoxPlus_{j \in R_i} F_{i,j} \beta_j \bigg)_{a_i}(1_{\bZ_2})\,.
\end{multline*}
where equality $(*)$ must hold for some coefficients $F_{i,j}$ with $j \in R_i$.
Setting $F_{i,j}=0$ for all $1\leq j\leq M$ with $j \not\in R_i$ leads to the associated matrix $F$ for $f^*$. In practise, the matrix $F$ is obtained by reducing a matrix 
\[
\Big( 
\Ima(d_{n+1}) \Big\vert \Cycle_n(K_*) \Big\vert f\Cycle_n(L_*)
\Big)\,,
\]
where $f\Cycle_n(L_*)$ denotes the set of embedded cycles. More precisely, tracking which columns from the second block were added to the columns from the third block leads to $F$. Also, notice that here we can reduce the embedded cycles only by the image and cycles that are smaller or equal to their filtration values. At the end of the reduction process, the block $f\Cycle_n(L_*)$ should be reduced to the zero matrix. 

\section{Quotient Computation Example}
\label{app:example-quotient}
In this appendix, we compute a quotient of persistence modules explicitly by following the algorithm presented in~\cite{Torras2023}. In particular, we compute the quotient that appears on the extension problem from Example~\ref{ex:extension_matrix}. First, consider the basis $\cE_1$ of $E^2_{1,0}$ given by the two generators from Exmample~\ref{ex:extension_matrix}, where we have $\beta_1\sim [0.59, 3.63)$ and $\beta_2\sim [0.60, 2.71)$. Here we write real numbers with a precision of two digits since in this example this was enough to distinguish most filtration values up to a small tolerance (for some values we keep an additional digit). On the other hand, we consider a basis sorted in the endpoint order $\cE_0=\{\gamma_i\}_{i=1}^8$ for the term $E^2_{0,1}$ where the associated intervals are 
\begin{align*}
    \gamma_1 & \sim  [0.68, 0.86) 
    & \hspace{1cm} 
    \gamma_4 & \sim  [2.71, 3.38) 
    \\
    \gamma_2 & \sim  [0.79, 1.006) 
    & \hspace{1cm} 
    \gamma_5 & \sim  [2.71, 12.21)
    & \hspace{1cm} 
    \gamma_7 & \sim  [1.63, 12.51) 
    \\
    \gamma_3 & \sim  [2.58, 2.77) 
    & \hspace{1cm} 
    \gamma_6 & \sim  [1.007, 12.35)
    & \hspace{1cm} 
    \gamma_8 & \sim  [3.63, 12.52)  & . 
\end{align*}
Next, we consider the persistence modules $\rB E^2_{0,1}$ and $\rD E^2_{0,1}$ introduced in Subsection~\ref{sub:extension-problem}, together with their respective barcode bases $\rB \cE^2_{0,1}=\{\beta_1^\rB, \beta_2^\rB\}$ and  $\rD \cE^2_{0,1}=\{\beta_1^\rD, \beta_2^\rD\}$. Thus, we have the associated intervals $\beta_1^\rB\sim [0.59, \infty)$, $\beta_2^\rB\sim [0.60, \infty)$ and $\beta_1^\rD\sim [3.63, \infty)$, $\beta_2^\rD\sim [2.71, \infty)$. Further, we order the generators from $\rB \cE^2_{0,1}$ and $\rD \cE^2_{0,1}$ by their birth values, so that $\beta^\rB_1 < \beta^\rB_2$ while $\beta^\rD_2 < \beta^\rD_1$. Recall that our aim is to compute the quotient of $E^2_{1,0}\oplus \rB E^2_{0,1}$ modulo $\textrm{ext}_1(\rD E^2_{1,0})$. We start by considering the matrix 
\begin{equation}\label{eq:block-matrix-quotient}
\left(
\begin{array}{c|cc}
\EXT_1(\cE_0, \cE_1) & 0 & P(\cE^2_{0,1})\\
P(\cE^2_{1,0}) & \Id_{\vert\cE^2_{1,0}\vert} & 0
\end{array}
\right)
\end{equation}
where $\EXT_1(\cE_0, \cE_1)$ is the extension matrix from Subsection~\ref{sub:extension-problem} and where the matrices $P(\cE^2_{1,0})$ and $P(\cE^2_{0,1})$ are the permutations from the standard to the endpoint order on the bases $\cE^2_{1,0}$ and $\cE^2_{0,1}$ respectively. In our particular case, matrix~(\ref{eq:block-matrix-quotient}) becomes
\begin{equation}\label{eq:quotient-block-matrix-particular}
\left( 
\begin{array}{r|cc|cc|cc|cc|ccc|c}
& \beta^\rD_2 & \beta^\rD_1 & \beta^\rB_1 & \beta^\rB_2 & 
\gamma_1 & \gamma_2 & \gamma_6 & \gamma_7 & \gamma_3 & 
\gamma_5 & \gamma_4 & \gamma_8
\\
\hline
 \gamma_1 & & & & & \redCell 1 & 0 & & & & & & \\
 \gamma_2 & & & & & 0 & \redCell 1 & & & & & & \\
 \hline 
 \gamma_3 & 1 & 0 & & & & & & & \redCell 1 & 0 & 0 & \\ 
 \gamma_4 & 1 & 0 & & & & & & & 0 & 0 & \redCell 1 & \\
 \gamma_5 & 1 & 0 & & & & & & & 0 & \redCell 1 & 0 & \\
 \hline 
 \gamma_6 & 1 & 0 & & & & & \redCell 1 & 0 & & & & \\
 \gamma_7 & 1 & 0 & & & & & 0 & \redCell 1 & & & & \\
 \hline 
 \gamma_8 & 0 & 1 & & & & & & & & & & \redCell 1\\
 \hline 
 \beta^B_1 & 0 & \redCell 1 & \redCell 1 & 0 & & & & & & & & \\
 \beta^B_2 & \redCell 1 & 0 & 0 & \redCell 1 & & & & & & & & 
\end{array}
\right)
\end{equation}
where the empty blocks are trivial. For ease, we have highlighted the pivots from matrix~(\ref{eq:quotient-block-matrix-particular}) using red-coloured cells. Here notice that the rows and the columns after the first column block are indexed by the union $\cE^2_{0,1}\cup \rB\cE^2_{1,0}$. Further, notice that such union is ordered following the standard order on the columns and the endpoint order on the rows; if this did not hold, one needs to perform the necessary permutations on~(\ref{eq:quotient-block-matrix-particular}) until the generators indexing both rows and columns are ordered. Here we remark that the first column block is independent of the columns on its right and should always precede the following columns. 

Now we are ready to apply the \texttt{box\_gauss\_reduce}~algorithm from~\cite{Torras2023}. We omit some of the details, like the identity matrix which keeps track of column additions or the use of the lists \texttt{lpivots} and \texttt{lpiv\_idx}. 
First, we record the birth values of columns into a list that we keep updating
\[
\lbirths = [
2.71, 3.63, 0.58, 0.60,0.67, 0.79,
1.007, 1.63, 2.37, 2.71, 2.71, 3.63
]
\]
and also store a copy of this list $\texttt{lbirths\_orig}$ which we do not alter. Also, we record the death values of rows in a list
\[
\ldeaths = [
 0.86, 1.006, 2.77, 3.38, 12.21, 12.35, 12.51, 12.52, \infty, \infty
]\,.
\]
Now, the \texttt{box\_gauss\_reduce}~algorithm reduces the matrix~(\ref{eq:quotient-block-matrix-particular}) by sweeping from the bottom to the top row.

We start from the bottom row (index $10$), and notice that both columns $1$ and $4$ have pivots on this row. We add the first column to the forth column and update the fourth entry of $\lbirths$ accordingly so that $\lbirths[4]\gets \max(\lbirths[1], \lbirths[4])(= 3.63)$.
Then, we might turn some entries in the fourth column to zero\footnote{
    For all rows $j < 10$, we check whether the entry in position $(j,4)$ from the matrix is non-trivial. If it is nontrivial we check that $\lbirths[4]<\ldeaths[i]$; if this does not hold we substitute this entry by a zero.
}  and move to the next row
and repeat this procedure. 
Eventually we obtain the matrix
\begin{equation}\label{eq:matrix-quotient-reduced}
\left( 
\begin{array}{cc|cc|cc|cc|ccc|c}
 & & &         & \redCell 1 & 0 & & & & & & \\
 & & &         & 0 &\redCell 1 & & & & & & \\
 \hline 
 1 & 0 & 0 & 1 & & & 0 & 1 &\redCell  1 & 1 & \blCell 0 & \\ 
 1 & 0 & 0 & 1 & & & 0 & 1 & 0 & \redCell 1 & \blCell 0 & \\
 1 & 0 & 0 & 1 & & & 0 & \redCell 1 & 0 & \blCell 0 & 0 & \\
 \hline 
 1 & 0 & 0 & 1 & & & \redCell 1 & \blCell 0 & & & & \\
 1 & 0 & 0 &\redCell  1 & & & 0 & \blCell 0 & & & & \\
 \hline 
 0 & 1 & \redCell  1 & 0 & & & & & & & & \blCell  0\\
 \hline 
 0 & \redCell 1 & \blCell 0 & 0 & & & & & & & & \\
 \redCell 1 & 0 & 0 & \blCell  0 & & & & & & & & 
\end{array}
\right)
\end{equation}
where the red highlighted entries are the resulting column pivots and the blue entries are those which have been pivots at some point during the algorithm.
Also, the resulting list of births is
\[
\lbirths = [
2.71, 3.63, \textbf{3.63}, \textbf{2.71}, 0.67, 0.79,
1.007, \textbf{2.71}, 2.37, \textbf{2.71}, 2.71, 3.63
]
\]
where we have highlighted in boldface those values which have been updated during the procedure.

Now we can read out the barcode of $\PH_1(\AXn)$ depicted in Subfigure~\ref{subfig:PH1-barcode-extension} from our computations. For this we range over the columns of indexes $j\geq 3$ from matrix~(\ref{eq:matrix-quotient-reduced}). There are two options for the $j^{\rm th}$-column. If the $j^{\rm th}$-column is non-trivial, we read the index of its pivot $i$ and obtain the interval $[\texttt{lbirths\_orig}[j], \ldeaths[i])$. Otherwise, if the $j^{\rm th}$-column is zero, we obtain the interval $[\texttt{lbirths\_orig}[j], \lbirths[j])$. Denoting by $I_j$ the barcode that we obtain for column $j$, the resulting intervals are 
\begin{align*}
    & 
    I_3 = [0.58, 12.52) 
    & & 
    I_6 = [0.79, 1.006)
    & & 
    I_9 = [2.37, 2.77)
    \\ & 
    I_4 = [0.60, 12.51) 
    & & 
    I_7 = [1.007, 12.35)
    & & 
    I_{10} = [2.71, 3.38)
    \\ & 
    I_5 = [0.67, 0.86) 
    & & 
    I_8 = [1.63, 12.21)
\end{align*}
and the trivial intervals $I_{11} = [2.71, 2.71)=\emptyset$ and $I_{12} = [3.63, 3.63)=\emptyset$. These are the intervals from the barcode depicted in Subfigure~\ref{subfig:PH1-barcode-extension}.
\vspace{1cm}
\end{document}